\documentclass[12pt]{amsart}
\usepackage[utf8]{inputenc}
\usepackage{amsmath,amsthm,amssymb}
\usepackage{xcolor}
\newcommand{\correction}[1]{{#1}}
\DeclareRobustCommand{\bluecorrection}[1]{{#1}}

\usepackage[british]{babel}

\usepackage[shortlabels]{enumitem}
\usepackage{hyperref}
\usepackage{mathtools}
\usepackage{graphicx}
\usepackage{tikz}
\usetikzlibrary{arrows.meta,decorations.pathreplacing}
\tikzset{>=Latex,every picture/.style={line cap=round,line join=round},
	point/.style={circle,fill=black,inner sep=1.35pt},
	lab/.style={font=\small},boxlabel/.style={font=\small,align=center},
	blueaxis/.style={inkblue,line width=1pt},
	greenaxis/.style={inkgreen,line width=1pt,dashed},
	orangeaxis/.style={inkorange,line width=1pt,dash dot}}
    
\hypersetup{colorlinks=true,linkcolor=black,citecolor=black,
	urlcolor=black,pdftitle={Three-ends Furstenberg incidences and approximate delta-rectangles},
	pdfauthor={Ciprian Demeter and William O'Regan},
	pdfsubject={Three-ends Furstenberg incidences and approximate delta-rectangle counting}}
\definecolor{inkblue}{RGB}{35,79,121}
\definecolor{inkorange}{RGB}{174,94,24}
\definecolor{inkgreen}{RGB}{35,114,93}
\definecolor{softblue}{RGB}{231,240,247}
\definecolor{softgreen}{RGB}{233,244,239}
\definecolor{softorange}{RGB}{252,241,225}

\newtheorem{theorem}{Theorem}[section]
\newtheorem{conj}[theorem]{Conjecture}

\newtheorem{remark}[theorem]{Remark}
\newtheorem{lemma}[theorem]{Lemma}
\newtheorem{proposition}[theorem]{Proposition}
\newtheorem{prop}[theorem]{Proposition}

\theoremstyle{definition}

\newtheorem{definition}[theorem]{Definition}

\mathtoolsset{showonlyrefs}
\newcounter{stepctr}
\newif\ifinsteppage

\newenvironment{steps}{%
	\setcounter{stepctr}{0}%
	\insteppagefalse
}{%
	\par\bigskip
}

\newcommand{\step}[1]{%
	\ifinsteppage
	\par\bigskip
	\fi
	\insteppagetrue
	\refstepcounter{stepctr}%
	\noindent\textbf{Step \thestepctr:}\if\relax\detokenize{#1}\relax\else\space#1.\fi\par
}

\newcommand{\dem}{\de^{-O(\e)}}
\newcommand{\dep}{\de^{O(\e)}}
\newcommand{\demm}{\de^{-\e}}
\newcommand{\depp}{\de^\e}

\numberwithin{equation}{section}

\newcommand{\R}{\mathbb{R}}
\newcommand{\T}{\mathcal{T}}

\newcommand{\Z}{\mathbb{Z}}

\newcommand{\Th}{\Theta}

\newcommand{\cT}{\mathcal{T}}

\newcommand{\cD}{\mathcal{D}}

\newcommand{\N}{\mathbb{N}}

\newcommand{\bT}{\mathbf{T}}

\newcommand{\dy}[2]{\cD_{#1}(#2)}

\newcommand{\cC}{\mathcal{C}}

\newcommand{\cR}{\mathcal{R}}

\newcommand{\en}[2]{\mathrm{E}_{#1,#2}}

\newcommand{\de}{\delta}
\newcommand{\e}{\epsilon}

\newcommand{\sg}{\sigma}

\newcommand{\bbP}{\mathbb{P}}

\newcommand{\cns}[2]{N_{{#1}}({#2})}

\newcommand{\cn}[1]{N_\de({#1})}

\newcommand{\card}[1]{\lvert #1\rvert}

\newtheorem*{ack*}{Acknowledgment}

\author{Ciprian Demeter}
\author{William O'Regan}
\begin{document}
	\title{Three-ends Furstenberg estimates and approximate rectangles}

	\thanks{CD is partially supported by NSF grant DMS-2349828. WOR is supported in part by an NSERC Alliance grant administered by Pablo Shmerkin and Joshua Zahl.}
	\begin{abstract}
\bluecorrection{We prove a three-ends incidence estimate whose multiplicity exponent is sharp for one-dimensional Katz--Tao point sets.} At the heart of the argument lies rectangular non-concentration and its incidence theory introduced earlier by the authors.   As an application, we obtain a strong quantitative bound on the number of approximate rectangles with the vertices belonging to four separated camps, each retaining a non-concentration condition. 
	\end{abstract}
	\maketitle
\section{Introduction}
Much of the work in approximate incidence geometry has centred on the study of Furstenberg sets. A \emph{$(\de,s,t,C)$-Furstenberg set} is a collection $\T$ of $\de$-tubes in the plane, where there exists a $(\de,t,C)$-set of points $P$ in the plane so that through each point $p \in P$ there exists a $(\de,s,C)$-set $\T(p) \subset \T$ of tubes incident to $p.$ 

\bluecorrection{For planar point configurations we work in a fixed bounded region, normalised to $[0,1]^2$. Point sets are $\de$-separated, and tube families are $\de$-separated in line parameter space, as specified in Section 2. All distances between line directions are measured modulo $\pi$.} 

The following sharp lower bound up to exponent for the number of tubes was obtained by \cite{renwang} building on the work of \cite{orpshabc}. We invite the reader to view the two aforementioned papers, along with the references therein, to view the colourful history of the problem of finding the sharp lower bound which has generated a rich body of mathematics over the last thirty years or so.
\begin{theorem}[Orponen--Shmerkin; Ren--Wang]\label{thm.furst}
\correction{Let $0<s\leq1$ and $0<t<2$. For every $\eta>0$ there exists $\epsilon_0=\epsilon_0(s,t,\eta)>0$ such that the following holds whenever $0<\epsilon\leq\epsilon_0$ and $\de>0$ is sufficiently small.}

Let $\cT$ be a $(\de,s,t,\demm)$-Furstenberg set. Let $M$ be an integer with the property that $|\T(p)| \geq M$ for each $p \in P.$
Then 
\begin{equation}
    |\cT| \mathrel{\correction{\gtrsim_{\eta}}}\correction{\de^{\eta}} M \de^{-\frak{f}(s,t)},
\end{equation}
where $$\frak f(s,t) =\min\left\{\frac{s+t}{2},t, 1\right\}.$$
\end{theorem}
Some further research directions have centred upon weakening some of the assumptions in the statement of Theorem \ref{thm.furst}. In one direction, one may ask what happens if the assumption on $P$ is weakened. If $P$ is only $\de$-separated, then the best one may obtain is 
\begin{equation}\label{eq.trivial}
|\T| \mathrel{\correction{\gtrsim_{\e}}}\correction{\dep} M\sqrt{|P|};
\end{equation}
see Proposition \ref{prop.bilinear}. This bound may be seen to be sharp if one considers a collection $P$ maximally clustered, that is, being contained in a ball of radius $\sim \de |P|^{1/2}.$ Provided this is not the case, in some weakly quantitative way, then one may improve upon \eqref{eq.trivial}. This was established in \cite{orponenshmerkin2026mild} and we invite the reader to view their well-written and clear paper.

In a different direction, one may choose to keep the non-concentration assumption on $P,$ but weaken or alter the assumption on the $\T(p).$ One example is to ask that each $\T(p)$ satisfies a certain {\em quantitative two-ends} assumption. \correction{In the case $t=1$, for any fixed $s\in(0,1]$ and any $\eta>0$, taking the input loss parameter $\e$ sufficiently small in terms of $s,\eta$ gives}
$$|\cT|\mathrel{\correction{\gtrsim_{s,\eta}}}\correction{\de^{\eta}} M^{3/2}\de^{-1/2}.$$
This is due to \cite{WW24}, generalised to other values of $t$ in \cite{wangwufurst}.
When $t \neq 1$ one also needs another notion of non-concentration on the $\T(p).$ See Theorem \ref{theorem.twoends} for the full statement. The hypotheses of their statement are not formally weaker, although morally they may be considered that way.

The question still remains: What is the weakest condition one may place on the $\T(p)$ to improve upon \eqref{eq.trivial} at the level of exponents? We provide a condition in Definition \ref{def.three-ends} below which is rather weak but sufficient.

\subsection{Incidence theorems}
\begin{definition}\label{def.three-ends}
    		Let $\T$ be a collection of $\de$-tubes, $P$ a collection of $\de$-separated points, $M \geq 1$ an integer, and $c >0.$ We say that $(P,\T)$ is \textit{$(M,c)$-three-ends} if for every $p\in P$, there are subfamilies
		\[
		\T_1(p),\T_2(p),\T_3(p)\subset\{T\in\T:T\cap p\ne\varnothing\},
		\]
        so that
$$
\card{\T_i(p)}\ge M \text{ for all } i = 1,2,3,
$$
and
		\begin{equation}\label{eq:local-transversality}
			\angle(T,T')\ge c
			\quad\text{if }T\in\T_i(p),\ T'\in\T_j(p),\ i\ne j.
		\end{equation}
        Here $\angle(T,T')$ is the angle between the centre-lines of $T$ and $T'.$
\end{definition}
\correction{If instead of three families we only have two,  $(P,\cT)$ will be referred to as $(M,c)$-two-ends. This is a non-quantitative version of the quantitative two-ends structure investigated earlier in literature. }

Set $$\sigma = \min \{s,2-s\}.$$
	\begin{theorem}\label{thm.3ends}
		Fix $0<s<2,\e >0.$ Let $P$ be a $(\de,s)$-KT-set with $|P| \sim \de^{-s}.$ Suppose that $(P,\T)$ is $(M,\depp)$-three-ends. Then
		\[
		|\T| \gtrsim_{\e}\dep M^{1+\sigma/4} \sqrt{|P|}
		\]
	\end{theorem}
The gain from three vs. two families of tubes can be seen by comparing it with Proposition \ref{prop.bilinear}.

For $s=1$, the exponent $5/4$ of $M$ in Theorem \ref{thm.3ends} is sharp (as a universal exponent that works for all $M$). The following multiscale construction was communicated to us by Tuomas Orponen. It shows the sharpness of $5/4$ when $M\sim \delta^{-2/3}$.

\begin{figure}[htbp]
		\centering
		\begin{tikzpicture}[scale=.94]
			\fill[black!1] (-2.40,-3.00) rectangle (2.40,3.00);
			\draw[black,line width=.95pt] (-2.40,-3.00) rectangle (2.40,3.00);
			
			\draw[black,line width=.8pt] (-2.40,-1.00)--(2.40,-1.00);
			\draw[black,line width=.8pt] (-2.40, 1.00)--(2.40, 1.00);
			
			\foreach \yy in {2.02,0.02,-1.98}{
				\draw[red!82!black,line width=.72pt] (-2.34,{\yy-.13})--(2.34,{\yy+.06});
				\draw[red!82!black,line width=.72pt] (-2.34,{\yy+.05})--(2.34,{\yy+.03});
				\draw[red!82!black,line width=.72pt] (-2.34,{\yy+.24})--(2.34,{\yy+.12});
			}
			
			\foreach \x/\y in {-1.88/1.30,-1.57/1.61,-1.26/1.92,-.95/2.23,
				-1.88/-.68,-1.57/-.37,-1.26/-.06,-.95/.25,
				-1.88/-2.66,-1.57/-2.35,-1.26/-2.04,-.95/-1.73}{
				\node[draw=black,line width=.8pt,minimum size=.34cm,inner sep=0pt,rotate=45] at (\x,\y) {};
			}
			
			\foreach \x/\y in {-.40/1.24,-.09/1.55,.22/1.86,.53/2.17}{
				\node[draw=black,line width=.8pt,minimum size=.34cm,inner sep=0pt,rotate=45] at (\x,\y) {};
			}
			
			\draw[black,line width=.78pt] (-3.70,-.38)--(-.66,2.48);
			\draw[black,line width=.78pt] (-3.70,-.52)--(-.66,2.52);
			\draw[black,line width=.78pt] (-3.70,-.68)--(-.66,2.56);
			
			\draw[inkblue,line width=.9pt] (-1.18,4.02)--(1.20,-4.02);
			\draw[inkblue,line width=.9pt] (.08,4.02)--(.42,-4.02);
			\draw[inkblue,line width=.9pt] (1.34,4.02)--(-.94,-4.02);
			
			\draw[<->,line width=.75pt] (-2.40,3.48)--(2.40,3.48)
			node[midway,above,lab] {$\delta^{1/3}$};
			\draw[<->,line width=.75pt] (-2.78,1.00)--(-2.78,3.00)
			node[midway,left,lab] {$\delta^{2/3}$};
			\draw[<->,line width=.75pt] (2.78,-3.00)--(2.78,3.00)
			node[midway,right,lab] {$1$};
			\node[lab,anchor=west] at (-2.24,2.62) {$\delta\text{-points}$};
			
			\node[boxlabel,anchor=east,text width=2.25cm] at (-3.62,.70)
			{$\delta^{-2/3}$ black\\tubes\\in one bush};
			\draw[->,line width=.65pt] (-3.54,.50)--(-3.02,.14);
			
			\node[boxlabel,anchor=west] at (3.48,2.14)
			{$\delta^{-1}$ red tubes\\in each horizon-\\tal rectangle};
			\draw[->,line width=.65pt] (3.38,1.95)--(2.05,2.17);
			
			\node[boxlabel,anchor=west,text width=2.25cm] at (3.48,-2.12)
			{$\delta^{-4/3}$ tubes\\of the blue family};
			\draw[->,line width=.65pt] (3.38,-2.05)--(1.18,-2.55);
		\end{tikzpicture}
		\caption{\bluecorrection{The sharpness example, redrawn from the communicated sketch.}}\label{fig:orponen}
	\end{figure}
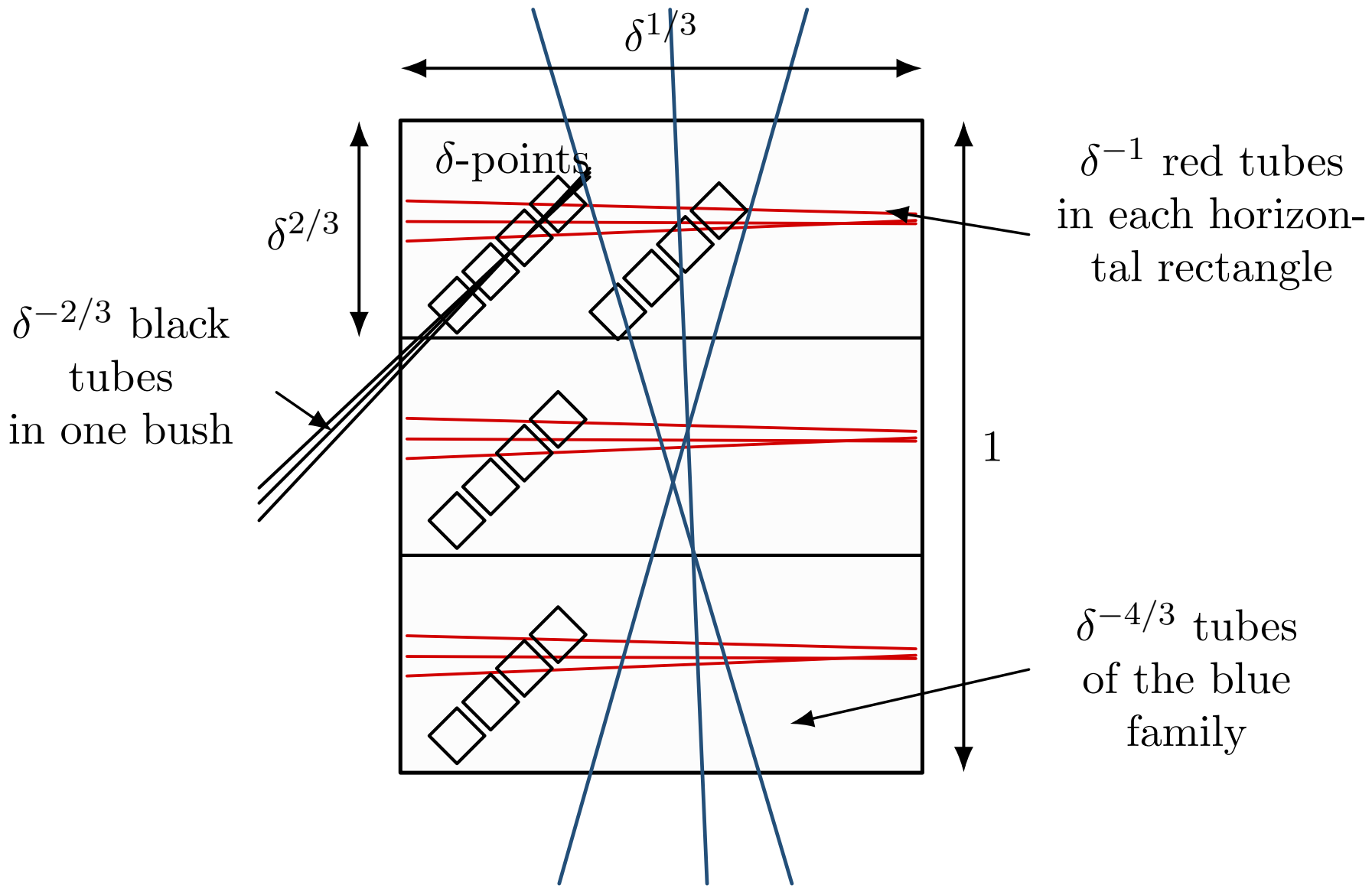

Set $\rho=\delta^{1/3}$. Its geometry is illustrated in Figure \ref{fig:orponen}. We pack a larger vertical block with thin horizontal rectangles. The three tube families are represented by lines: the red lines run more or less along the long directions of the horizontal rectangles (with the small angular variation permitted by their eccentricity), the blue lines run along the long direction of the larger vertical block, and the dark lines go through one of the clusters. The three direction families are quantitatively transverse.
	
	The scales are chosen so that
	\begin{equation}\label{eq:Orponen-counts}
		\card{P}\sim\delta^{-1},\qquad
		\card{\T_i}\sim \delta^{-4/3},
		\qquad
		\card{\T_i(p)}\sim M=\delta^{-2/3}
		\quad(p\in P)
	\end{equation}
	for each of the three colours. \bluecorrection{For an explicit configuration, take an integer $n\geq2$, put $\delta=n^{-3}$, $\rho=n^{-1}$, and
    \[
    P=\{(a\rho^2+b\delta,c\rho+b\delta):a,b,c\in\{0,\ldots,n-1\}\}.
    \]
    Rows are separated vertically by $\rho$, clusters within a row horizontally by $\rho^2$, and points within a cluster by $\sqrt2\delta$ along a $45^\circ$ line. The bound $|P\cap B(x,r)|\lesssim r/\delta$ follows by considering $\delta\leq r\leq\rho^2$, $\rho^2\leq r\leq\rho$, and $\rho\leq r\leq1$. For the red, black and blue families, respectively, take the tubes meeting $P$ from a $\delta$-net in line parameter space with direction windows of width $c_0\rho$ about $0,\pi/4,\pi/2$, where $c_0>0$ is a small fixed constant. Each point meets $\sim\rho/\delta=\delta^{-2/3}$ tubes of each colour. There are $\sim\delta^{-1}$ red tubes per row, $\sim\delta^{-2/3}$ black tubes per cluster, and $\sim\delta^{-4/3}$ blue tubes in total. These counts give \eqref{eq:Orponen-counts}.}
	
	Indeed,
	\[
	\delta^{-1/2}M^{5/4}
	=\delta^{-1/2}(\delta^{-2/3})^{5/4}
	=\delta^{-4/3},
	\]
	which is exactly the size of each tube family in \eqref{eq:Orponen-counts}. Thus the exponent $5/4$ cannot be improved in general without adding new hypotheses. 
    
    We state the following more general version where $P$ need not have the `correct' cardinality for its dimension, and may have a large Katz--Tao constant.
    \begin{proposition}\label{prop.main}
        		Fix $0<s<2,C \geq 1, \e >0.$ Let $M \geq 1$ be an integer. Let $P$ be a $(\de,s,C)$-KT-set. Suppose that $(P,\T)$ is $(M,\depp)$-three-ends. Then
		\[
		|\T| \gtrsim_{\e}C^{-3/4} M^{1+\sigma/4} \delta^{s/4+O(\e)}|P|^{3/4}.
		\]
    \end{proposition}
    We will prove Proposition \ref{prop.main} when $C \sim 1.$ The general case follows easily from an application of \cite[Lemma 2.7, Lemma 2.8]{demwangszem}. Indeed, after an application of both, one finds a $P' \subset P$ with $|P'| \gtrsim \dep |P|/C,$ which is a $(\de,s,O_\e(1))$-KT set. Applying Proposition \ref{prop.main} with $P'$ gives the desired result.

    Let us now show some applications, that have motivated our study of the three-ends phenomenon. 

    \subsection{Approximate rectangles}
    Let $S_1,S_2,S_3,S_4$ be disjoint sets in $[0,1]^2.$ We say that 
    $$(p_1,p_2,p_3,p_4) \in S_1 \times S_2 \times S_3 \times S_4$$
    is a \emph{$\de$-approximate rectangle} if 
    \begin{equation}\label{eq.para}
    |p_1-p_2+p_3-p_4| = O(\de),
    \end{equation}
    and
    \begin{equation}\label{eq.angle}
    |(p_2-p_1)\cdot (p_4-p_1)| = O(\de).
    \end{equation}
    Certainly this definition matches what one may expect a $\de$-approximate rectangle to be: \eqref{eq.para} tells us $(p_1,p_2,p_3,p_4)$ forms a `$\de$-parallelogram'; \eqref{eq.angle} tells us that the sides formed from $p_1$ make an angle of $\pi/2 \pm O(\de/(|p_2-p_1||p_4-p_1|))$, which is $\pi/2 \pm O(\de^{1-2\e})$ under the separation assumption below. 
    
    \bluecorrection{Using the mixed circle--tube variant of Proposition \ref{prop.main} proved in Section \ref{sect.rect}, we obtain the following result.} 
    \begin{theorem}\label{thm.rect}
        Let $\e >0.$ The following holds for all $\de > 0$ small enough.

        Let $S_1,S_2,S_3,S_4 \subset [0,1]^2$ be $(\delta,1)$-KT sets, each of size $\de^{-1},$ and mutually separated by a distance of $\depp.$ Let $\cR$ be the $\de$-approximate rectangles formed from $S_1,S_2,S_3,S_4.$ Then
        \begin{equation}
            |\cR| \lesssim_\e \dem\correction{\de^{-7/3}}.
        \end{equation}
    \end{theorem}
    We will prove the more general version below.
    \begin{prop}\label{prop.rect}
        Let $0 < s < 2$ and let $\e >0$ be sufficiently small. The following holds for all $\de > 0$ small enough.

        Let $S_1,S_2,S_3,S_4 \subset [0,1]^2$ be $(\delta,s)$-KT sets, each of size at most $N,$ and mutually separated by a distance of $\depp.$ Let $\cR$ be \bluecorrection{the collection of $\de$-approximate rectangles} formed from $S_1,S_2,S_3,S_4.$ Then
        \begin{equation}
        	\label{eq.fjvhuhvutyrtugu999t}
            |\cR| \lesssim_\e \de^{-\frac{s}{4-\sigma}-O(\e)}N^{\frac{9-3\sigma}{4-\sigma}}
        \end{equation}
        When $s=1$, this is
        $$\bluecorrection{|\cR|\lesssim_\e\delta^{-1/3-O(\e)}N^2.}$$
    \end{prop}
    The exact rectangle analogue is classical. Pach and Sharir proved that, for any fixed angle, $n$ planar points determine at most $O(n^2\log n)$ ordered triples realizing that angle~\cite{PS}; in particular, exact rectangles admit an essentially quadratic bound by looking at a right-angle corner. The thickened problem is substantially more delicate. An exact right-angle line becomes a $\delta$-tube, and the Katz-Tao condition  does not prevent many points from concentrating in a long thin rectangle inside such a tube. Thus the non-concentration can only be exploited in the forms permitted by available Furstenberg-type incidence theorems. The three-ends result is designed precisely to manufacture that stronger form of non-concentration. To take advantage of it, we encode rectangles as hyperedges, uniformize various local degrees and do counting relative to substructures of the hypergraph.

      We refer the reader to the last section for a broader discussion.

    \subsection{Energy estimates on the paraboloid}
    The $\delta$-approximate rectangle count admits an equivalent formulation. 
    Denote the paraboloid in $\R^3$ by
\begin{equation}
    \bbP^2 = \{(x,y,x^2+y^2) : |x|, |y| \leq 1\}. 
\end{equation}
\begin{theorem}
    Let $0 < s < 2, \e >0.$ The following holds for all $\de >0$ small enough. Let $P_1,P_2,P_3,P_4 \subset \bbP^2$ be $(\de,s)$-KT sets, mutually separated by a distance of $\depp.$
    Denoting 
    $$\en{4}{\de}(P_1,P_2,P_3,P_4) := |\{(p_1,p_2,p_3,p_4) \in P_1 \times P_2 \times P_3 \times P_4 : |p_1 -p_2 + p_3 - p_4| = O(\de)\},$$
    we have
    \begin{equation}\label{eq.energyresult}
        \en{4}{\de}(P_1,P_2,P_3,P_4) \lesssim_\e\dem \de^{-\frac{s(10-3\sg)}{4-\sigma}}.
    \end{equation}
\end{theorem}
We sketch how this immediately follows from Proposition \ref{prop.rect}. 
There is a correspondence between $\delta$-approximate quadruples on the paraboloid and $\delta$-approximate  rectangles. 
Select $$(p_1,p_2,p_3,p_4) \in P_1 \times P_2 \times P_3 \times P_4,$$
which satisfies 
\begin{equation}\label{eq.energy}
|p_1 - p_2 + p_3 - p_4| = O(\de).
\end{equation}
Let $s_i$ be the image of $p_i$ under the canonical projection to $\R^2;$ clearly for $(s_1,s_2,s_3,s_4)$ equation \eqref{eq.para} is satisfied. Equation \eqref{eq.angle} is satisfied since 
\begin{align}
 |(s_2-s_1)\cdot (s_4-s_1)|
    &= \frac12\bigl||s_1|^2-|s_2|^2+|s_3|^2-|s_4|^2\bigr|+O(\de)\\
    &= O(\de),
\end{align}
where we have used that \eqref{eq.energy} is satisfied for the third coordinate. Since this canonical projection is bi-Lipschitz with constant $\sim 1,$ the $S_i$ are $(\de,s)$\correction{-KT sets}, and still mutually separated by a distance of $\sim \depp,$ we may apply Proposition \ref{prop.rect} to obtain \eqref{eq.energyresult}.

The reader is referred to Section \ref{sect.fourth-final} for a discussion of energies. We will see that when $s\le 1$, our exponent \eqref{eq.energyresult} for the separated energy outperforms the optimal exponent for the linear energy. Also, the sharpness of the upper bound \eqref{eq.energyresult} deteriorates for $s>1$, we only record it for completeness. In fact, the sharp exponent when $s\ge 3/2$ follows via a more elementary approach, as explained in the last section.

The case $s=1$ seems most appropriate for testing the strength of our approach. The exponent $7/3$ uses our three-ends result which, as observed earlier, has a sharp exponent $5/4$ for $M$. Further tiny improvements of $7/3$ are possible by using Theorem \ref{theorem.twoends} in the extreme case when all the  three-ends-based  inequalities are saturated. We will not pursue this improvement here, as it  involves unappealing numerics. Getting the optimal exponent 2 would be a significant target.   

The interest in the energy formulation goes back to  a closely related $\delta$-discretized energy problem on the parabola
\[
\mathbb P^1=\{(t,t^2):|t|\le 1\}.
\]
Let $0<s\le \frac12$, and let $P\subset\mathbb P^1$ be $\delta$-separated and satisfy the $s$-dimensional nonconcentration condition
\[
\#(P\cap B(x,r))\lesssim (r/\delta)^s,
\qquad \delta\le r\le1.
\]
Its $\delta$-discretized six-point energy is
\[
\en{6}{\delta}(P)
:=\#\bigl\{(p_1,\ldots,p_6)\in P^6:
|p_1+p_2+p_3-p_4-p_5-p_6|\lesssim\delta\bigr\}.
\]
The sharp conjectural estimate
\begin{equation}\label{eq:parabola-six-energy-conj}
	 \en{6}{\delta}(P)\lesssim_\varepsilon \delta^{-3s-\varepsilon}
\end{equation}
was first posed in this form by Orponen~\cite{O1}. His work gave a power-saving improvement over the elementary bounds, which exposed connections to deep incidence geometry. The authors later proved
\begin{equation}\label{eq:parabola-six-energy-DO}
 \en{6}{\delta}(P)\lesssim_\varepsilon \delta^{-\frac72 s-\varepsilon}.
\end{equation}
See~\cite{doprod} and the references therein for the intervening developments.
    \subsection{Projection theorems}
    We write $S^1_+ := \{(\cos\theta,\sin\theta): 0 \leq \theta \leq \pi/2\}$ for the first-quadrant arc.
	\begin{theorem}\label{thm.proj}
		Fix $0<s<2,\e >0.$ Let $P\correction{\subset[0,1]^2}$ be a $(\de,s)$-KT-set with $|P| \sim \de^{-s}.$ Let $\Theta \subset S^1_+$ \bluecorrection{be $\de$-separated} and suppose there exist $x_1,x_2,x_3 \in \Theta$ so that 
        $$|x_i-x_j| \geq \depp \text{ for } i \neq j.$$
        Then there exists $x \in \Theta$ so that 
		$$
		\cn{\pi_x(P)} \gtrsim_{\e}\dep |\Th|^{\sigma/4} \sqrt{|P|}.
		$$
	\end{theorem}
    \correction{A qualitative version of Theorem \ref{thm.proj}, assuming additionally $|\Theta|\geq\de^{-\tau}$ for a fixed $\tau>0$ and taking $\e$ sufficiently small in terms of $s,\tau$, yields a gain $\de^{-\eta}$ for some $\eta=\eta(s,\tau)>0$. This was obtained by Ren and Shmerkin and shared with us via private communication.} They did not have access to the recalled Theorem \ref{thm.rectangle} \cite{doprod} and instead elected to opt for an application of \cite[Theorem 1.7]{shmerkin2023nonlinear}, see also \cite[Theorem 1.2]{orponenshmerkin2026mild}. Their proof strategy has some similarities to ours.

We  prove the following more general version.
    	\begin{prop}\label{prop.proj}
		Fix $0<s<2,\e >0.$ Let $P\correction{\subset[0,1]^2}$ be a $(\de,s,C)$-KT-set. Let $\Theta \subset S^1_+$ \bluecorrection{be $\de$-separated} and suppose that there exist $x_1,x_2,x_3 \in \Theta$ so that 
        $$|x_i-x_j| \geq \depp \text{ for } i \neq j.$$
        Then there exists $x \in \Theta$ so that 
		\begin{equation}\label{eq.propconcl}
		\cn{\pi_x(P)} \gtrsim_{\e}C^{-3/4}\dep |\Th|^{\sigma/4} \de^{s/4}|P|^{3/4}.
		\end{equation}
	\end{prop}
    Again, we prove this when $C \sim 1,$ with the reduction to the general case following from the short argument presented below the statement of Proposition \ref{prop.main}. \correction{When $C\sim1$, \eqref{eq.propconcl} improves the bilinear baseline $\dep\sqrt{|P|}$ by a fixed power of $\de^{-1}$, provided that}
    \begin{equation}
        \correction{|\Th|^{\sigma}|P|\geq\de^{-s-\kappa}}
    \end{equation}
    \correction{for some fixed $\kappa>0$, with $\e$ sufficiently small in terms of $\kappa$.}
    We remark Proposition \ref{prop.proj} does not follow from Proposition \ref{prop.main}.

    \subsection*{Acknowledgements}
    We are indebted to Kevin Ren and Pablo Shmerkin for sharing a proof of a non-quantitative version of Theorem \ref{thm.proj}. We thank Tuomas Orponen for sharing a sharp example for Theorem \ref{thm.3ends}, as well as bringing the matter of Theorem \ref{thm.3ends} to the attention of WOR. CD thanks Shukun Wu for  pointing out to him the train-track example relevant to two-ends. WOR thanks Pablo Shmerkin and Joshua Zahl for support.
    
    This material is based upon work supported by the Swedish Research Council under grant number 2021-06594 whilst WOR was in residence at Institut Mittag-Leffler in Djursholm, Sweden during the autumn of 2026. Thanks are given to the organisers, administrators and the participants.
    \subsection*{Statement on the usage of artificial intelligence}
    The proofs of Theorem \ref{thm.3ends} (Proposition \ref{prop.main}) and Theorem \ref{thm.proj} (Proposition \ref{prop.proj}) were found by the authors. The approach of Theorem \ref{thm.rect} (Proposition \ref{prop.rect}) was found by the authors. A proof of Proposition \ref{prop.rect} using Proposition \ref{prop.main} was attempted. ChatGPT-6 was then employed to check the proofs. It found a substantial error in the original proof of Proposition \ref{prop.rect}, and kindly suggested that it may be fixed by modifying the statement of Proposition \ref{prop.main} to Proposition \ref{prop.circle}, which employs a proof which is essentially only different on a cosmetic level. More details may be found in Section \ref{sect.rect}. The manuscript was written by humans although we have used AI in an assistive role, which did include writing some short and prescriptive parts. For example, all the text in the next section preceding the statement of Proposition \ref{prop.bilinear} and part of the content of Section \ref{sect.fourth-final} was generated by AI, before being edited by humans. The pictures were AI generated.
    
\section{Preliminaries}

We use $|P|$ for the cardinality of a finite set or collection. Write $N_r(P)$ for the least number of radius $r$ balls needed to cover $P$. For dyadic $r$, let $\dy r P$ be the collection of half-open dyadic $r$-cubes meeting $P$; thus $|\dy r P|\sim N_r(P)$. All point sets are $\de$-separated, so $N_\de(P)\sim |P|$. We freely identify a point with a ball or square of size $\de$ about it, interpreting incidences up to an $O(\de)$ error. For $\theta\in S^1$, write $\pi_\theta(x)=x\cdot\theta$.

The notation $A\lesssim B$ means $A\leq CB$, with $C$ independent of $\de$ and $\e$; dependence on fixed dimension parameters is allowed. Subscripts indicate additional dependence. We write $A\gtrsim B$ for $B\lesssim A$, and $A\sim B$ when both hold. Also, $f=O(g)$ means $|f|\lesssim g$. In exponents, $O(\e)$ denotes a non-negative quantity bounded by $C\e$, with $C$ independent of $\de,\e$ and allowed to change between occurrences. Thus $\dem$ and $\dep$ denote losses by small powers of $\de$. 

We recall some non-concentration conditions. For $0\leq s\leq d$ and $C\geq1$, a $\de$-separated set $P\subset\R^d$ is a \emph{$(\de,s,C)$-set} if
\[
    |P\cap B(x,r)|\leq Cr^s|P|\qquad \text{for } x\in\R^d \text{ and } r\geq\de.
\]
It is a \emph{$(\de,s,C)$-KT set}, or Katz--Tao set, if instead
\[
    |P\cap B(x,r)|\leq C(r/\de)^s\qquad\text{for } x\in\R^d \text{ and }r\geq\de.
\]
We also recall the following concept introduced in  \cite[Definition  1.4]{doprod}. 
\begin{definition}
A planar $P$ is a \emph{$(\de,s,C)$-RKT set}, or rectangular Katz--Tao set, if
\[
    |P\cap R|\leq C(\sqrt{rw}/\de)^s
\]
for every rectangle $R$ with side lengths $\de\leq r\leq w$, in any orientation.
\end{definition}
 We omit $C$ when $C\sim1$. For sets which are not discrete, the same definitions use $N_\de$ in place of cardinality.

A $\de$-tube is the $\de$-neighbourhood of a line, restricted to a fixed bounded region. Tube families are $\de$-separated in line parameter space. Locally, after rotation, we identify the line $y=ax+b$ with $(a,b)$; non-concentration conditions on tubes refer to these parameters. For tubes incident to one point, this is equivalent, up to constant factors, to the corresponding condition on their directions.

We call $P$ uniform at prescribed dyadic scales if, at each scale $\rho$, the quantities $|P\cap Q|$ are comparable independently of $Q\in\dy\rho P$. We use scales separated by a fixed factor $2^{-T}$\bluecorrection{, apart from a possibly shorter terminal block}, with $T=T(\e)$ sufficiently large. A uniform subset of size $\gtrsim\de^\e|P|$ exists by \cite[Lemma 2.6]{doprod}. At intermediate dyadic scales the comparison constants may depend on $\e$; in particular $|P\cap Q|\sim_\e |P|/N_\rho(P)$ for occupied $Q$.

    \begin{proposition}\label{prop.bilinear}
        Let $M \geq 1, c > 0.$ Suppose that $(P,\T)$ is $(M,c)$-two-ends. Then 
        \begin{equation}
            |\cT| \gtrsim c M |P|^{1/2}.
        \end{equation}
    \end{proposition}
    The proof is elementary, and may be found in many locations, such as \cite[Proposition 2.9]{orpyi}. Note that we do not need any non-concentration on $P.$ 
    \begin{proof}[Proof of Proposition \ref{prop.bilinear}]
        We have
        \begin{equation}
            M^2|P| \leq |\{(p,T,S)\in P \times \T \times \T : p\cap T\ne\varnothing,\ p\cap S\ne\varnothing \text{ and } \angle(T,S) \geq c\}| \lesssim c^{-2}|\T|^2.
        \end{equation}
        The lower bound is immediate; the upper bound follows since for any pair $(T,S)$ with the conditions stipulated above, the admissible $p$ will live in a ball of radius $\sim \de/c,$ which may trivially contain $\lesssim c^{-2}$ points from $P.$ The result follows immediately. 
    \end{proof}
    For $c \sim 1,$ Proposition \ref{prop.bilinear} is sharp, even when $P$ has some non-concentration.  An example would be the well-known train-track example. One should view \cite[Construction 1; figure 3]{fu2021incidence}. Let $P$ be \bluecorrection{a maximal $\de$-separated subset of a union of} $\sim \de^{-1/2} \ $ $\de \times \de^{1/2}$ rectangles, the short end situated on the $x$-axis, each separated by distance $\de^{1/2}.$ Let $1 \leq M \leq \de^{-1/2}.$ For each $p \in P,$ consider a family of $M$ tubes $T_1(p)$ with maximal concentration, pointing vertically, and a family of $M$ tubes $T_2(p),$ with maximal concentration, pointing horizontally. Typically, a tube is incident to $\de^{-1/2}$ points. Double counting the incidences, we arrive at
    $$M\de^{-1} \sim |\T|\de^{-1/2},$$
    showing that Proposition \ref{prop.bilinear} is sharp.

    We recall the dual form of \cite[Theorem 4.4]{oreganwuyi2026} in the case where the transversal family of curves consists of lines. This is a standard corollary of  \cite[Theorem 1.5]{doprod}, which was proved with $C,K\sim 1$.
\begin{theorem}\label{thm.rectangle}
		Let $0 < s < 2, \e >0.$ The following holds for all $\de >0$ small enough.
        
        Consider a collection $\T$ of $\delta$-tubes and a collection $P$ of $\delta$-points in $[0,1]^2$. Assume that  $P$ is a $(\delta,s,C)$-RKT set. Assume that for each $p\in P$, a subfamily $Y(p)\subset\{T\in\T:\;p\cap T\not=\emptyset\}$ is a $(\delta,\sigma,K)$-KT set, $\sigma=\min(s,2-s)$. Suppose $|Y(p)| \geq M$ for all $p \in P.$
Then
		\begin{equation}
			\label{ orjifu54it0-96-py=56-56}
			|\T| \geq \frac{\dep}{C\sqrt{K}} M^{3/2}\de^{s/2} |P|.
		\end{equation}
\end{theorem}
We also require the `quantitative two-ends Furstenberg' theorem of Wang and Wu, \cite[Theorem 1.5]{wangwufurst}. We state it in the dual form.
\begin{theorem}
    \label{theorem.twoends}
    Let $P$ be a $(\delta,s,C)$-KT set of
    $\delta$-balls in $[0,1]^2,$
    for some $s\in(0,2)$, and let $\e >0$.
    Write $\sigma=\min(s,2-s)$.

    For each $p\in P$, let $\T(p)$ be a family of $\delta$-tubes meeting $p$.
    Assume that each $\T(p)$ is a $(\delta,\sigma,K)$-KT
    set and also a $(\delta,\epsilon,\de^{-\epsilon^{3}})$-set,
    and that $|\T(p)|\geq M$.

    Write $\T=\bigcup_{p\in P}\T(p).$ Then
    $$
        |\T|
        \gtrsim_{\epsilon}
        \frac{\dep}{ C\sqrt K}
        M^{3/2}\delta^{s/2}|P|.
    $$
\end{theorem}
We recall the ubiquitous `two-ends reduction'. For example, one may find a similar version at \cite[Lemma 2.4.6]{TaoEpsilonII}. We give a rendition in the form we require.
\begin{lemma}\label{lem.twoends}
 Let $\e >0.$ The following holds for all $\de >0$ small enough. Let $\Theta \subset S^1$ be $\de$-separated. There exist a scale $\rho \in [\de,1]$ and an arc $I$ of length $\rho$ so that $|\Theta \cap I| \mathrel{\correction{\gtrsim}} \de^{\e}|\Theta|,$ and $S_I(\Theta \cap I)$ is a $(\de/\rho, \e, \correction{O(1)})$-set. Here $S_I$ is \correction{the normalised arclength parametrisation of $I$ onto $[0,1]$}.
\end{lemma}
\correction{\begin{proof}
Choose $I$, with $\rho=|I|\in[\de,1]$, maximising $|\Theta\cap I|/\rho^\e$. Covering $S^1$ by $O(1)$ unit-length arcs gives $|\Theta\cap I|\gtrsim\rho^\e|\Theta|\geq\de^\e|\Theta|$. For any subarc $J\subset I$ of length at least $\de$, maximality gives $|\Theta\cap J|\leq(|J|/\rho)^\e|\Theta\cap I|$. Rescaling proves the stated Frostman bound.
\end{proof}}

\section{Proof of Proposition \ref{prop.main}}
{
\noindent\textit{Proof outline.} 
The argument exploits bilinear transversality through the three-ends hypothesis. For every rectangle, two of the three separated direction families are transverse both to each other and to its long axis. This permits bilinear incidence counting inside rectangles of arbitrary orientation and eccentricity. Combined with pruning at all scales, these estimates produce a quantified rectangular Katz–Tao condition, to which we apply our earlier incidence theorem, more precisely its version with explicit dependence on the RKT constant, Theorem \ref{thm.rectangle}. 

\subsection{Pre-processing}
Throughout it will be implicitly assumed that $\de,\e$ are taken small enough to ensure validity.

Firstly, it is convenient to reduce to a more global condition regarding the directions of the tubes lying through each $p \in P.$

\begin{lemma}\label{lem.intervals}
    There exist $I_1,I_2,I_3 \subset S^1,$ with $d(I_j,I_k) \gtrsim \depp$ for $j\ne k$, and $P_0 \subset P$ with $|P_0| \geq \dep|P|$ so that for each $p \in P_0,$ there are $$\T_1(p),\T_2(p), \T_3(p) \subset \{T\in\T:T\cap p\ne\varnothing\}$$ with the slopes of $\T_1(p),\T_2(p), \T_3(p)$ belonging to $I_1,I_2,I_3$ respectively. Furthermore, $|\T_i(p)| \geq M'$ for each $i = 1,2,3$, for an integer $M'\gtrsim\depp M$.
    \begin{proof}
        Tile $S^1$ with arcs of length at most $\depp/10$ and comparable to $\depp$:
        $$\{I_1,\dots, I_N\},$$
        where $N \sim \demm.$ For each $p$, by the pigeonhole principle, restrict each of its three tube families to an arc containing at least $M'=\lceil M/N\rceil$ of its tubes. We have
        \begin{align}
            |P| = \sum_{\substack{i,j,k\text{ pairwise distinct}}}\sum_{p \in P} 1_{\sigma(\T_1(p)) \times \sigma(\T_2(p)) \times \sigma (\T_3(p)) \subset I_i \times I_j \times I_k},
        \end{align}
        where $\sigma(\ \cdot \ )$ denotes slope. We may therefore find pairwise distinct $i,j,k$ and a $P_0 \subset P$ with $|P_0| \geq \dep |P|$ so that for all $p \in P_0$ we have that
        $$\sigma(\T_1(p)) \times \sigma(\T_2(p)) \times \sigma (\T_3(p)) \subset I_i \times I_j \times I_k.$$
        Note that since $(P,\T)$ is $(M,\depp)$-three-ends we automatically acquire the necessary separation for the intervals.
    \end{proof}
\end{lemma}
From now on we replace $P$ with $P_0$ and $M$ with $M'$. We may also now write $\T = \T_1 \cup \T_2 \cup \T_3$ where the slopes from each $\T_j$ belong to $I_j$ for $j = 1,2,3.$ We now replace $P$ by a uniform subset of size at least $\depp |P|.$

\bluecorrection{The next lemma controls the products of local tube counts at every coarse scale.}

For $i = 1,2,3$ and $Q \in \dy \rho P$ write
$$\T_i(Q) = \bigcup_{p \in P \cap Q}\T_i(p).$$  
\begin{lemma}\label{lem.rectangularpruninglol}
There exists $P_1 \subset P$ with $|P_1| \geq \dep|P|$ so that the following holds for all $\de \leq \rho \leq 1.$ 

Let $Q \in \dy \rho {P_1}.$ Then 
\begin{equation}
     |\T_i(Q)||\T_j(Q)| \leq \dem \frac{|\T|^2}{\cns \rho {P_1}},
\end{equation}
for all $i, j \in \{1,2,3\}$ with $i < j.$
\begin{proof}
   We begin with $\{i,j\} = \{1,2\}.$ Obviously one may repeat this for other pairs by re-labelling.
   
   Without loss of generality we may suppose $\de=2^{-n}$ is dyadic. \bluecorrection{Use the same sufficiently large integer block length $T=T(\e)\geq1/\e$ as in the uniformisation, put $m=\lceil n/T\rceil$, and set $\rho_j=2^{-\min\{jT,n\}}$ for $0\leq j\leq m$.} Let $\rho = \rho_j$ for some $1 \leq j \leq m.$ Let $Q \in \dy \rho P.$ We have
   \begin{align}
   \sum_{Q \in \dy \rho P}|\T_1(Q)||\T_2(Q)| \leq \de^{-10\e}|\T_1||\T_2|.
   \end{align}
   Indeed, fix $(T_1,T_2)\in\T_1\times\T_2$. Since $\rho\geq\de$,
every $\rho$-square meeting both tubes lies in the intersection of
the $O(\rho)$-neighbourhoods of their centre-lines. Their angular
separation places this intersection inside a ball of radius
$O(\rho\de^{-\e})$, which meets at most $\de^{-10\e}$ such squares, provided that $\de$ is small enough.
Thus each tube pair contributes at most $\de^{-10\e}$ times to
the sum. Summing over all pairs gives the above.

It is straightforward to see that for all but at most $\de^\e\cns\rho P$ of the $Q \in \dy \rho P$ we have 
\begin{equation}\label{eq.good}
    |\T_1(Q)||\T_2(Q)| \leq \de^{-11\e} \frac{|\T|^2}{\cns \rho P}.
\end{equation}
If not then this implies
\begin{equation}
      \de^{-10\e}|\cT|^2 < \de^{-10\e}|\T|^2,
\end{equation}
which is nonsense.

Apply this argument at each $1 \leq j \leq m,$ starting from $j=1$, for all three pairs, using the current surviving set and its covering number in each application of \eqref{eq.good}. Since earlier deletions remove whole squares at this or coarser scales, uniformity shows that each application removes at most an $O_\e(\de^\e)$ fraction of the remaining points, which is at most a tenth for $\de$ small enough. Thus at least half the points survive each scale. The resulting set $P_1 \subset P$, with $P$ here denoting the original set, satisfies
\begin{equation}\label{eq.boundP1}
|P_1| \geq 2^{-m}|P|.
\end{equation}

Now for $\de < \rho < 1,$ choose $0 \leq j \leq m-1$ so that
$$\rho_{j+1}\leq \rho \leq \rho_j.$$
For $Q \in \dy \rho {P_1}$ let $Q' \in \dy {\rho_j} {P_1}$ be such that $Q \subset Q'.$ Every intermediate surviving set contains $P_1$, so, using the trivial bound when $j=0$,
\begin{equation}\label{eq.conclusion}
    |\T_1(Q)||\T_2(Q)| \leq |\T_1(Q')||\T_2(Q')| \leq \de^{-11\e}\frac{|\T|^2}{\cns {\rho_j} {P_1}} \leq 4^T\de^{-11\e} \frac{|\T|^2}{\cns \rho {P_1}}.
\end{equation}
\bluecorrection{Since $T\geq1/\e$, we have $m\leq\e n+1$. Thus \eqref{eq.boundP1} gives}
$$\bluecorrection{|P_1|\geq\tfrac12\de^\e|P|\geq\dep|P|,}$$
and \eqref{eq.conclusion} may be bounded by
$$\de^{-100\e} \frac{|\T|^2}{\cns \rho {P_1}},$$
as required. 
\end{proof}
\end{lemma}
We again abuse notation and set $P = P_1.$
\subsection{Proof of Proposition \ref{prop.main}}
The aim is to apply Theorem \ref{thm.rectangle} directly. The parameters relating to the rectangularity of $P$ are stated and proved in the proposition below. 
\begin{proposition}\label{prop.RKT}
    The collection $P$ is a $(\de,s,C)$-RKT set, where
$$C \lesssim \dem \frac{|\T|}{M|P|^{1/2}}.$$
\end{proposition}
\begin{proof}
    Let $\de \leq r \leq w,$ with $r < \depp,$ and let $R$ be a rectangle with shorter side $r$ and longer side $w.$ (If $w \geq r \geq \depp$ then the KT bound of $(w/\de)^s$ will suffice here.) The aim is to count
    $$|P \cap R|.$$
    (Cover $R$ by $O(1)$ dyadic $w$-squares; the argument below is applied in each such square $Q$, and the resulting bounds are summed.)
    We may cover the occupied part of $R$ by $N \lesssim w/r$ $r$-squares from $\dy r P;$ call them $S_i$ for $1 \leq i \leq N.$ We then have, using the $(\de,s)$-KT condition on $P$ that
    \begin{align}
        |P \cap R| &\leq \sum_{i=1}^N |P \cap S_i|\\
        &\leq \max_{1 \leq i \leq N}|P \cap S_i|^{1/2} \sum_{i=1}^N |P \cap S_i|^{1/2}\\
        &\lesssim (r/\de)^{s/2}  \sum_{i=1}^N |P \cap S_i|^{1/2}.
    \end{align}
    
    The aim is to now bound 
    \begin{equation}\label{eq.boundthis}
    \sum_{i=1}^N |P \cap S_i|^{1/2}.
        \end{equation}
 By the three-ends assumption, each $p \in P$ intersects three mutually $\depp$-separated families of tubes. By Lemma \ref{lem.intervals}, we may suppose, without loss of generality, that the tubes from $\T_1$ and $\T_2$ make an angle $\gtrsim\depp$ from the centre-line of $R.$ Applying the bilinear argument from Proposition \ref{prop.bilinear} we have
 \begin{equation}
     M^2|P \cap S_i|\leq \dem |\T_1(S_i)||\T_2(S_i)|.
 \end{equation}
 Since $r< \depp,$ any $T \in \T_j$ may only intersect $O(\demm)$ of the squares $S_i$, for $j=1,2.$ We may therefore bound
\begin{align}
     \sum_{i=1}^N|P\cap S_i|^{1/2}
     &\lesssim \frac{\dem}{M} \sum_{i=1}^N|\cT_1(S_i)|^{1/2}|\cT_2(S_i)|^{1/2}\\
     &\leq \frac{\dem}{M}\sqrt{\sum_{i=1}^N|\cT_1(S_i)|\sum_{i=1}^N|\cT_2(S_i)|}\\
     &\leq \frac{\dem}{M} \sqrt{|\T_1(Q)||\T_2(Q)|}\\
     &\leq \frac{\dem}{M\sqrt {\cns w P}}|\T|.
\end{align}
Here we have used Cauchy--Schwarz and Lemma \ref{lem.rectangularpruninglol}. Since \correction{$P$ is a $(\de,s)$-KT set} it follows readily that
\begin{equation}
    \cns w P \geq \dep \de^s |P| w^{-s}.
\end{equation}
We may further bound \eqref{eq.boundthis} by
\begin{equation}
    \dem (w/\de)^{s/2}\frac{|\T|}{M|P|^{1/2}}.
\end{equation}
So we finally bound
\begin{equation}
    |P \cap R| \leq C (\sqrt{rw}/\de)^s
\end{equation}
 where
 $$C \lesssim \dem \frac{|\T|}{M|P|^{1/2}}.$$
 Incidentally, up to a factor $\dem$, $C$ is bounded by the ratio of improvement of $|\T|$ over what one obtains from the bilinear argument alone. 
\end{proof}
\begin{figure}[htbp]
		\centering
		\begin{tikzpicture}[scale=.95]
			\foreach \shift/\rot/\captiontext in {-3.1/0/{long axis near $\T_1$},3.1/60/{long axis near $\T_2$}}{
				\begin{scope}[xshift=\shift cm]
					\begin{scope}[rotate=\rot]
						\filldraw[fill=black!5,draw=black!65] (-2.05,-.35) rectangle (2.05,.35);
						\draw[densely dotted] (-2.25,0)--(2.25,0);
					\end{scope}
					\foreach \a in {-4,0,4}{
						\draw[inkorange!75,line width=.6pt] ({-1.7*cos(\a)},{-1.7*sin(\a)})--({1.7*cos(\a)},{1.7*sin(\a)});
						\draw[inkblue!75,line width=.6pt] ({-1.7*cos(60+\a)},{-1.7*sin(60+\a)})--({1.7*cos(60+\a)},{1.7*sin(60+\a)});
						\draw[inkgreen!75,line width=.6pt] ({-1.7*cos(120+\a)},{-1.7*sin(120+\a)})--({1.7*cos(120+\a)},{1.7*sin(120+\a)});
					}
					\fill (0,0) circle (1.8pt);
					\node[lab,inkorange] at (1.94,.18) {$\T_1$};
					\node[lab,inkblue] at (1.13,1.60) {$\T_2$};
					\node[lab,inkgreen] at (-1.13,1.60) {$\T_3$};
					\node[boxlabel] at (0,-2.13) {\captiontext};
				\end{scope}
			}
			\node[boxlabel] at (-3.1,-2.63) {use families $2,3$};
			\node[boxlabel] at (3.1,-2.63) {use families $1,3$};
		\end{tikzpicture}
		\caption{For each rectangle, the useful pair of families may change, but it is fixed for that rectangle and works at every one of its points. At most one of three separated directions can be nearly parallel to the long axis.}\label{fig:three-cones}
	\end{figure}
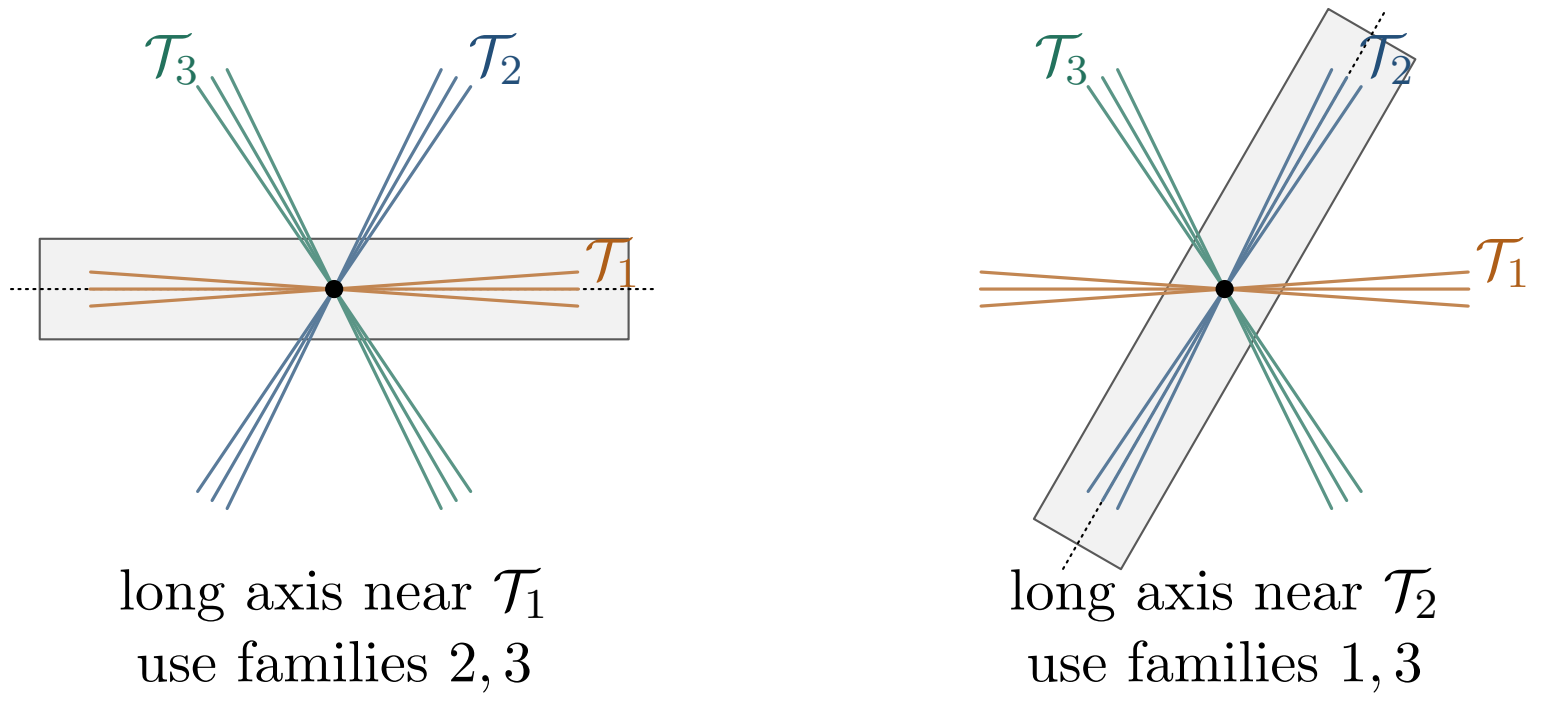
    	\begin{figure}[htbp]
		\centering
		\begin{tikzpicture}[scale=.97]
			\fill[softblue] (3.3,0) rectangle (4.4,1.1);
			\draw[thick] (0,0) rectangle (9.9,1.1);
			\foreach \x in {1.1,2.2,3.3,4.4,5.5,6.6,7.7,8.8}{\draw[black!35] (\x,0)--(\x,1.1);}
			\foreach \x/\y in {.24/.24,.76/.76,1.31/.45,1.90/.86,2.48/.72,2.95/.22,3.50/.34,3.90/.78,4.13/.49,4.72/.72,5.14/.28,5.78/.34,6.15/.81,6.83/.61,7.37/.20,8.03/.79,8.53/.36,9.01/.57,9.56/.28}{
				\fill (\x,\y) rectangle ++(.055,.055);
			}
			\draw[blueaxis] (3.10,-.6)--(4.18,1.9);
			\draw[greenaxis] (4.60,-.6)--(3.22,1.9);
			\node[lab,inkblue] at (4.65,1.85) {$T_1$};
			\node[lab,inkgreen] at (2.78,1.85) {$T_2$};
			\draw[<->] (0,-.37)--(9.9,-.37) node[midway,below,lab] {$w$};
			\draw[<->] (-.35,0)--(-.35,1.1) node[midway,left,lab] {$r$};
			\draw[decorate,decoration={brace,amplitude=4pt}] (3.3,1.28)--(4.4,1.28)
			node[midway,above=5pt,lab] {$r$};
			\node[boxlabel] at (8.7,2.05) {$|P \cap S_i|^{1/2} \lesssim \dem |\T_1(S_i)|^{1/2}|\T_2(S_i)|^{1/2}/M$};
			\draw[->,black!65] (6.0,1.80) to[out=180,in=60] (4.50,1.05);
		\end{tikzpicture}
		\caption{The one-rectangle count. A useful transverse tube sees only $\dem$ consecutive $r$-squares, even when $w/r$ is large. See Proposition \ref{prop.RKT}.}\label{fig:slabs}
	\end{figure}
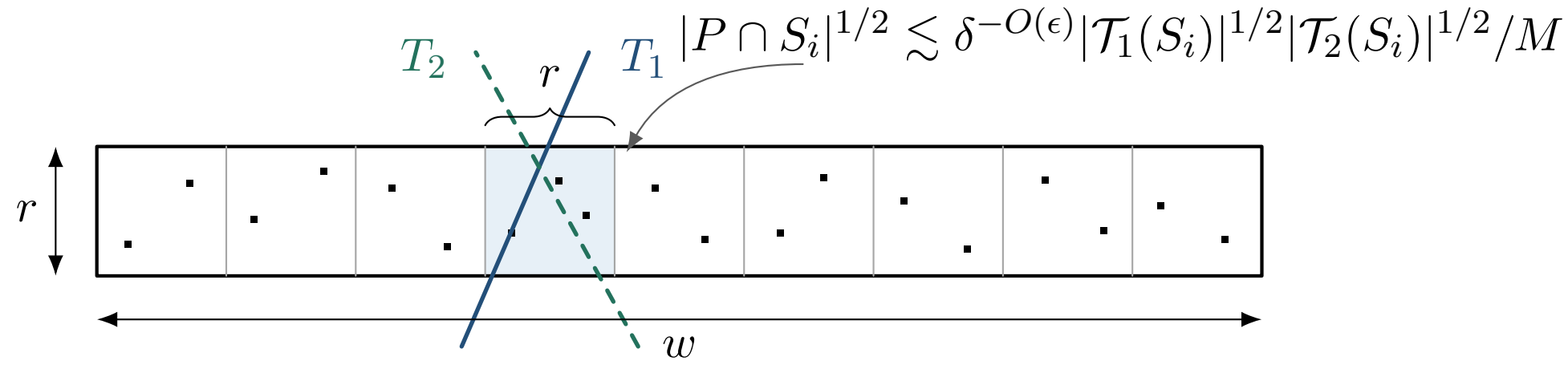
\begin{proof}[Proof of Proposition \ref{prop.main}]
    Choose $M$ tubes through each $p$ and denote their family by $\T(p)$. To apply Theorem \ref{thm.rectangle} we require that $\T(p)$ is a $(\de,\sigma,K)$-KT-set for some $K \geq 1.$ Trivially, letting $\bT$ be an $r \times 1$ tube we have 
    \begin{equation}
        |\{T\in\T(p):T\subset\bT\}| \lesssim \min\{(r/\de),M\} \lesssim M^{1-\sigma}(r/\de)^\sigma.
    \end{equation}
    Therefore applying Theorem \ref{thm.rectangle} with $K\lesssim M^{1-\sigma}$ we obtain that 
    \begin{equation}
        |\T|^2 \geq \dep M^{2 + \sigma/2}\de^{s/2}|P|^{3/2}, 
    \end{equation}
    which concludes the result.
\end{proof}

\section{Proof of Proposition \ref{prop.rect}}\label{sect.rect}
To prove Proposition \ref{prop.rect}, we require a variant of Proposition \ref{prop.main} where one bundle of tubes is replaced by a family of circles. We must also pay attention to a larger range of parameters, leading to a more involved statement.
\begin{proposition}\label{prop.circle}
    Let $0 < s < 2$ and let $M,L$ be positive integers. The following holds for all $\de, \e > 0$ small enough.

    Let $P \subset [0,1]^2$ be a $(\de,s)$-KT-set. Let $I_1,I_2,J$ be arcs in $S^1,$ mutually separated by a distance of $\depp$\correction{ in the angular metric modulo $\pi$}. Let $\T_1,\T_2$ be families of $\de$-tubes whose slopes lie in $I_1$ and $I_2$ respectively. Let $\cC$ be a family of $\de$-circles, with diameter between $\de^\e$ and $1.$ \correction{For each $C\in\cC$, write $\Gamma_C$ for its supporting circle and set
    \[
    \Gamma_C^J=\{q\in\Gamma_C:\operatorname{dir}(T_q\Gamma_C)\in J\},
    \]
    with directions read modulo $\pi$. This is a union of at most two connected circle arcs, and every tangent along each arc has direction in $J$. Only the $O(\de)$-neighbourhood of $\Gamma_C^J$ is used for circle incidences; $|\cC|$ still counts the underlying circle labels.}

    Suppose that for each $p \in P$ the following holds:
    \begin{enumerate}
        \item There is a subfamily $\T_1(p) \subset \T_1,$ all of whose members are incident to $p,$ with $|\T_1(p)| \geq M.$
        \item There is a subfamily $\T_2(p) \subset \T_2,$ all of whose members are incident to $p,$ with $|\T_2(p)| \geq M.$
        \item \correction{There is a subfamily $\cC(p)\subset\cC$ with $|\cC(p)|\geq L$ such that $\operatorname{dist}(p,\Gamma_C^J)=O(\de)$ for every $C\in\cC(p)$.}
    \end{enumerate}
    Then
    \begin{equation}
        \Sigma\times \min\{|\T_1|,|\T_2|\}^2\gtrsim_{\e}\dep M^{2+\sigma}\de^s|P|^3,
    \end{equation}
    where
    $$\Sigma=\frac{|\T_1||\T_2|}{M^2}+\frac{(|\T_1|+|\T_2|)|\cC|}{ML}.$$
\end{proposition}
	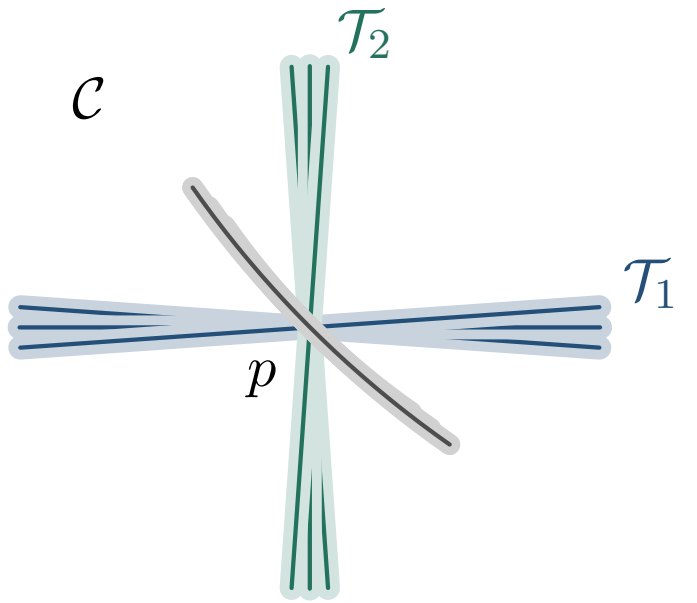
\begin{figure}[htbp]
		\centering
		\begin{tikzpicture}[scale=.96]
			\node[point] at (0,0) {};
			\node[lab,below left] at (-.08,-.08) {$p$};
			\foreach \a in {-4,0,4}{
				\draw[inkblue!25,line width=4.5pt] ({-2.0*cos(\a)},{-2.0*sin(\a)})--({2.0*cos(\a)},{2.0*sin(\a)});
				\draw[inkblue,line width=.8pt] ({-2.0*cos(\a)},{-2.0*sin(\a)})--({2.0*cos(\a)},{2.0*sin(\a)});
				\draw[inkgreen!20,line width=4.5pt] ({-1.8*sin(\a)},{-1.8*cos(\a)})--({1.8*sin(\a)},{1.8*cos(\a)});
				\draw[inkgreen,line width=.8pt] ({-1.8*sin(\a)},{-1.8*cos(\a)})--({1.8*sin(\a)},{1.8*cos(\a)});
			}
			\foreach \Rr in {3.7,4.4,5.1}{
				\draw[black!18,line width=4pt] plot[domain=215:235,samples=60]
				({\Rr+1.414*\Rr*cos(\x)},{\Rr+1.414*\Rr*sin(\x)});
				\draw[black!70,line width=.8pt] plot[domain=215:235,samples=60]
				({\Rr+1.414*\Rr*cos(\x)},{\Rr+1.414*\Rr*sin(\x)});
			}
			\node[lab,inkblue] at (2.35,.30) {$\T_1$};
			\node[lab,inkgreen] at (.38,2.02) {$\T_2$};
			\node[lab] at (-1.52,1.58) {$\cC$};
			\node[boxlabel,anchor=west,text width=5.9cm] at (3.2,.15)
			{$M$ selected tubes from each straight family;\\[5pt]
				$L$ selected circle arcs.};
		\end{tikzpicture}
		\caption{The mixed three-ends configuration. Blue and green distinguish only the two straight tube families; the circle arcs are drawn in grey.}\label{fig:mixed}
	\end{figure}
The proof of Proposition \ref{prop.circle} is essentially identical to the proof of Proposition \ref{prop.main}. We follow the proof there and emphasise the differences. \correction{All circle incidences are restricted as above, and for any square $Q$ we write $\cC(Q)=\bigcup_{p\in P\cap Q}\cC(p)$.} We now replace $P$ by a uniform subset of size at least $\depp|P|$.

The only difference in Lemma \ref{lem.rectangularpruninglol} is that a pair may now consist of a tube and a circle. \correction{We use the following elementary observation. If every tangent to a connected arc $\gamma$ makes an angle at least $c$ with a line $\ell$, then any strip of width $a$ parallel to $\ell$ cuts out an arc of length $O(a/c)$. Indeed, for a unit-speed parametrisation $\gamma(t)$ and a unit normal $n$ to $\ell$, the derivative $n\cdot\gamma'(t)$ has constant sign and absolute value $\gtrsim c$.

For a fixed tube--circle pair, every $\rho$-square counted by $|\T_j(Q)||\cC(Q)|$ lies within $O(\rho)$ of a retained arc in an $O(\rho)$-strip about the tube's centre-line. Applying the observation to the at most two arcs, with $c\gtrsim\depp$, places these squares in $O(1)$ balls of radius $O(\rho/c)$. Thus each pair is counted at most $O(c^{-2})=O(\de^{-2\e})$ times. Together with the tube--tube estimate, this gives, for all $\de\leq\rho\leq1$,}
\[
\sum_{Q\in\dy\rho P}\left(\frac{|\T_1(Q)||\T_2(Q)|}{M^2}
+\frac{(|\T_1(Q)|+|\T_2(Q)|)|\cC(Q)|}{ML}\right)\lesssim\dem\Sigma.
\]
The uniformisation step is identical: we obtain $P_1\subset P$ with $|P_1|\geq\dep|P|$ such that
\begin{equation}\label{eq.circlepruning}
\frac{|\T_1(Q)||\T_2(Q)|}{M^2}+\frac{(|\T_1(Q)|+|\T_2(Q)|)|\cC(Q)|}{ML}\leq\dem\frac{\Sigma}{\cns\rho{\bluecorrection{P_1}}}
\end{equation}
for all $Q\in\dy\rho{P_1}$. We again abuse notation and set $P=P_1$.
We now prove a version of Proposition \ref{prop.RKT} which only differs superficially. We sketch the difference.
\begin{proposition}\label{prop.circlerect}
The collection $P$ is a $(\de,s,C)$-RKT set, where
$$C\lesssim\dem\sqrt{\frac{\Sigma}{|P|}}.$$
\end{proposition}
\begin{proof}

    Let $\de \leq r \leq w,$ with $r < \depp,$ and let $R$ be a rectangle with shorter side $r$ and longer side $w.$
    \correction{The case $r\geq\depp$ follows from the KT bound: the tube--tube double count gives $\Sigma\gtrsim\de^{2\e}|P|$, so the additional factors are absorbed into $\dem$.}
    (Cover $R$ by $O(1)$ dyadic $w$-squares; the argument below is applied in each such square $Q$, and the resulting bounds are summed.)
    We may cover the occupied part of $R$ by $N \lesssim w/r$ $r$-squares from $\dy r P;$ call them $S_i$ for $1 \leq i \leq N.$ As before we have
    \begin{align}
        |P \cap R| &\lesssim (r/\de)^{s/2}\sum_{i=1}^N |P \cap S_i|^{1/2}.
    \end{align}
    As before we bound 
    \begin{equation}\label{eq.boundthiscircle}
    \sum_{i=1}^N |P \cap S_i|^{1/2}.
        \end{equation}

 By the assumptions, each $p \in P$ intersects three mutually $\depp$-separated families of tubes or circle arcs. We record the case where the tubes from $\T_1$ and the tangents from $\cC$ make an angle $\gtrsim\depp$ from the centre-line of $R.$ The case of two tube families is the same as in Proposition \ref{prop.RKT}. \correction{The strip estimate with $a=O(\de)$ shows that a fixed tube--circle pair has at most $O(\de^{-2\e})$ common selected $\de$-points. The bilinear double count from Proposition \ref{prop.bilinear} therefore gives}
 \begin{equation}
     ML|P \cap S_i|\leq \dem |\T_1(S_i)||\cC(S_i)|.
 \end{equation}
 \correction{The same strip estimate with $a=O(r)$, now parallel to the long axis of $R$, shows that each circle contributes to at most $O(\demm)$ of the squares $S_i$. The same bound holds for each tube from $\T_1$.} We may bound this sum by
\begin{align}
     \sum|P\cap S_i|^{1/2}
     &\lesssim \frac{\dem}{\sqrt{ML}} \sum|\T_1(S_i)|^{1/2}|\cC(S_i)|^{1/2}\\
     &\leq \frac{\dem}{\sqrt{ML}}\sqrt{\sum|\T_1(S_i)|\sum|\cC(S_i)|}\\
     &\leq \frac{\dem}{\sqrt{ML}} \sqrt{|\T_1(Q)||\cC(Q)|}\\
     &\leq \frac{\dem\sqrt{\Sigma}}{\sqrt {\cns w P}}.
\end{align}
Here we have used Cauchy--Schwarz and \eqref{eq.circlepruning}. Since $P$ is $(\de,s)$-KT it follows that
\begin{equation}
    \cns w P \geq \dep \de^s |P| w^{-s}.
\end{equation}
We may further bound \eqref{eq.boundthiscircle} by
\begin{equation}
    \dem (w/\de)^{s/2}\sqrt{\frac{\Sigma}{|P|}}.
\end{equation}
We may therefore finally bound
\begin{equation}
    |P \cap R| \leq C (\sqrt{rw}/\de)^s
\end{equation}
 where
 $$C \lesssim \dem \sqrt{\frac{\Sigma}{|P|}}.$$
 In conclusion $P$ is a $(\de,s,C)$-RKT-set.
\end{proof}
\begin{proof}[Proof of Proposition \ref{prop.circle}]
    For each $j=1,2$, choose $M$ tubes from $\T_j$ through each $p$ and denote their family by $\T(p)$. As before $\T(p)$ is a $(\de,\sigma,K)$-KT-set for $K \lesssim M^{1-\sigma}.$
    Therefore applying Theorem \ref{thm.rectangle} to the $\T_j$ we obtain
    \begin{equation}
        \Sigma\times \min\{|\T_1|,|\T_2|\}^2 \gtrsim_{\e} \dep M^{2+\sigma}\de^s|P|^3, 
    \end{equation}
    which concludes the result.
\end{proof}
\subsection{Proof of Proposition \ref{prop.rect}}
\begin{steps}
    \step{Setup}
    Let $\cR$ be the rectangles formed from $S_1,S_2,S_3,S_4.$ \bluecorrection{We use cyclic indices modulo $4$.} 

    By \eqref{eq.para} and \eqref{eq.angle}, all four vertices lie within $O(\de^{1-\e})$ of the circle \correction{having the line segment $[p_2,p_4]$ as diameter}. We first enlarge the working scale to $\de^{1-2\e}$. Partitioning into separated sets at this scale and applying \cite[Lemmas 2.8 and 2.9]{demwangszem}, we partition each $S_i$ into $\dem$ Katz--Tao sets at the new scale. Pigeonholing retains a $\dep$ fraction of the rectangles. Relabel the new scale as $\de$; the separation remains at least $\de^{O(\e)}$, and returning to the original scale only costs a factor $\dem$. A fixed rescaling ensures that all circle diameters are at most $1$.

    Let $\T_{i}$ be the $\de$-tubes which connect a point from $S_{i}$ to $S_{i+1}.$ Each $R \in \cR$ may be encoded by four tubes, one from each of the $\T_i,$ where counting is done modulo $4.$ Let $\cC$ consist of the $\de$-circles obtained by discretising the circles \correction{having the line segment $[p_2,p_4]$ as diameter}. At this scale the two diagonal circles differ by $O(\de)$ in their centres and radii, so either diagonal determines at most $O(1)$ possible encodings. Choose these encodings once and for all.

    At each corner, the two side directions and \correction{the tangent to $\Gamma_C$ at a point within $O(\de)$ of that corner} are separated by at least $\de^{O(\e)}$, since both adjacent sides have length at least $\de^{O(\e)}$. \correction{The circle diameters are at least $\de^{O(\e)}$, so the $O(\de)$ discretisation error changes these tangent directions by only $O(\de^{1-O(\e)})$, preserving the separation.} Pigeonholing their directions into sufficiently short arcs, we retain a further $\dep$ fraction of the rectangles so that these directions lie in three fixed, mutually separated arcs at each corner. \correction{Denote the circle-direction arc at $S_i$ by $J_i$ and use $\Gamma_C^{J_i}$ for all circle incidences with $S_i$. We keep the full-circle label $C$ in the hypergraph; at each corner it supplies at most two connected arcs.} Thus the angular and diameter assumptions of Proposition \ref{prop.circle} hold with $\e$ replaced by a constant multiple of $\e$.

    Each $R \in \cR$ may be encoded uniquely by the 9-tuple
    $$h \in H$$
    where 
    $$h = (p_1,p_{2},p_{3},p_{4}, T_1,T_2, T_3, T_4, C) $$
    \correction{Here $H$ is the set of these chosen encodings, so}
    $$H \mathrel{\correction{\subset}} S_1 \times S_2 \times S_3 \times S_4 \times \T_{1} \times \T_{2} \times \T_{3} \times \T_{4} \times \cC.$$
    The $p$ are the vertices of $R;$ the $T$ are the tubes enclosing the rectangle; $C$ is a circle containing the $p.$
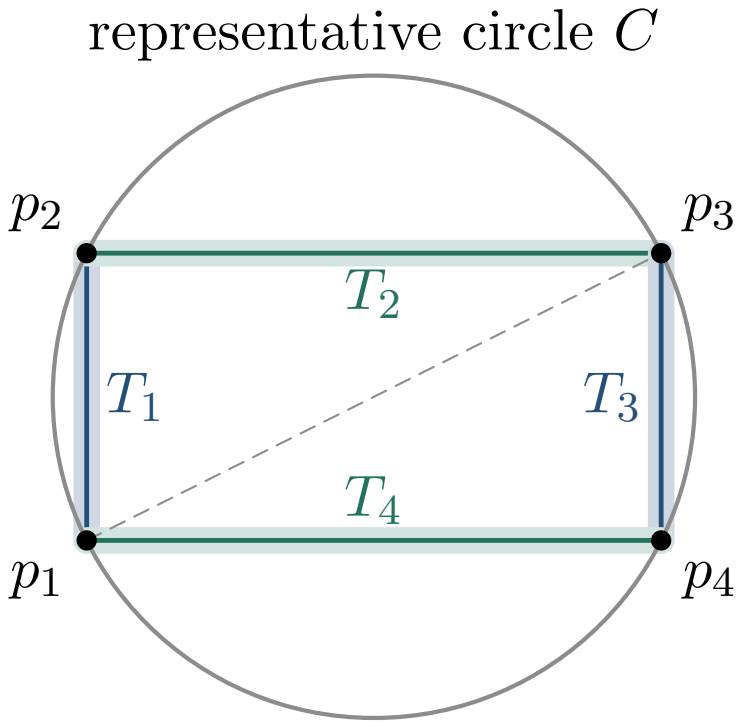
    	    	\begin{figure}[htbp]
		\centering
		\begin{tikzpicture}[scale=.95]
			\coordinate (p1) at (-2,-1); \coordinate (p2) at (-2,1);
			\coordinate (p3) at (2,1); \coordinate (p4) at (2,-1);
			\draw[inkblue!22,line width=5pt] (p1)--(p2);
			\draw[inkblue,line width=.9pt] (p1)--(p2);
			\draw[inkgreen!20,line width=5pt] (p2)--(p3);
			\draw[inkgreen,line width=.9pt] (p2)--(p3);
			\draw[inkblue!22,line width=5pt] (p3)--(p4);
			\draw[inkblue,line width=.9pt] (p3)--(p4);
			\draw[inkgreen!20,line width=5pt] (p4)--(p1);
			\draw[inkgreen,line width=.9pt] (p4)--(p1);
			\draw[black!45,line width=.8pt] (0,0) circle ({sqrt(5)});
			\draw[black!45,densely dashed] (p1)--(p3);
			\foreach \p/\lab/\pos in {p1/$p_1$/below left,p2/$p_2$/above left,p3/$p_3$/above right,p4/$p_4$/below right}{
				\node[point] at (\p) {}; \node[lab,\pos] at (\p) {\lab};}
			\node[lab,inkblue] at (-1.66,0) {$T_{1}$};
			\node[lab,inkgreen] at (0,.72) {$T_{2}$};
			\node[lab,inkblue] at (1.66,0) {$T_{3}$};
			\node[lab,inkgreen] at (0,-.72) {$T_{4}$};
			\node[lab] at (0,2.52) {representative circle $C$};
			\node[boxlabel,anchor=west,text width=6.3cm] at (3.3,.1)
			{One hyperedge is the nine-entry vector\\[5pt]
				$(p_1,p_2,p_3,p_4,T_{1},T_{2},T_{3},T_{4},C)$.\\[7pt]
				The points are black; blue and green distinguish only the two side-tube directions.};
		\end{tikzpicture}
        \end{figure}
    \step{Uniformisation}
    \bluecorrection{We use a variant of the hypergraph refinement argument of Wang--Zahl \cite[Lemma 8.1]{wang2022sticky}, combined with dyadic pigeonholing, to make the relevant nonzero extension degrees comparable.}
    We regard $H$ as the hypergraph whose edges are the tuples encoding rectangles. All incidences below are required to occur in a retained edge of $H$.

    For $p \in S_i$ let $M_i(p),M'_i(p),L_i(p)$ count the incident tubes from $\T_i,\T_{i-1}$ and circles from $\cC$, respectively. For $T \in \T_i$ let $r_i(T),r'_i(T)$ count the incident points from $S_{i+1},S_i$, respectively. For $C \in \cC$ let $t_i(C)$ count the incident points from $S_i$.

    Consider the nine kinds of vertices and the twelve kinds of pairs
    $$(p_i,T_i),\qquad (p_i,T_{i-1}),\qquad (p_i,C),\qquad i=1,2,3,4.$$
    For each such vertex or pair $u$, let $D(u)$ count the edges containing it. For each kind in turn, dyadically pigeonhole the current values of $D(u)$, retain the corresponding edges, and record the selected lower bound for each surviving $u$. This retains $\gtrsim (\log(1/\de))^{-21}|H|$ edges.

    Now repeatedly remove all edges containing any $u$ whose nonzero remaining degree is less than $\de^\e$ times its recorded lower bound. For each kind, the sum of these recorded lower bounds is at most the original $|H|$. Thus at most $21\de^\e|H|$ edges are removed in total, leaving a $\dep$ fraction of the original edges. The nonzero degrees of each kind now differ by at most $2\de^{-\e}$.

    Since, for example,
    $$D(p_i)=\sum_{T_i:\,(p_i,T_i)\text{ occurs in a retained edge}}D(p_i,T_i),$$
    the incidence counts are also uniform up to factors $\dem$. We may therefore choose
    $$(M_i,M'_i,r_i,r'_i,L_i,t_i)_{i=1}^4\in\N^{24}$$
    so that, for every retained edge and each $i$,
    \[
    \begin{aligned}
    \dep(M_i,M'_i,r_i,r'_i,L_i,t_i)
    &\leq (M_i(p_i),M'_i(p_i),r_i(T_i),r'_i(T_i),L_i(p_i),t_i(C))\\
    &\leq\dem(M_i,M'_i,r_i,r'_i,L_i,t_i).
    \end{aligned}
    \]
    Relabel the retained rectangles and edges as $\cR$ and $H$, and remove unused points, tubes and circles. Each $|S_i|$ remains at most $N$. We suppress these losses from now on, absorbing them into small powers of $\de$.

    For each $p \in S_i$, every retained tube from $\T_i(p)$ has a partner from $\T_{i-1}(p)$ making an angle $\pi/2+O(\de^{1-O(\e)})$, and conversely. A fixed tube through $p$ has at most $\dem$ such partners. Hence $\dep M_i\lesssim M'_i\lesssim\dem M_i$.

	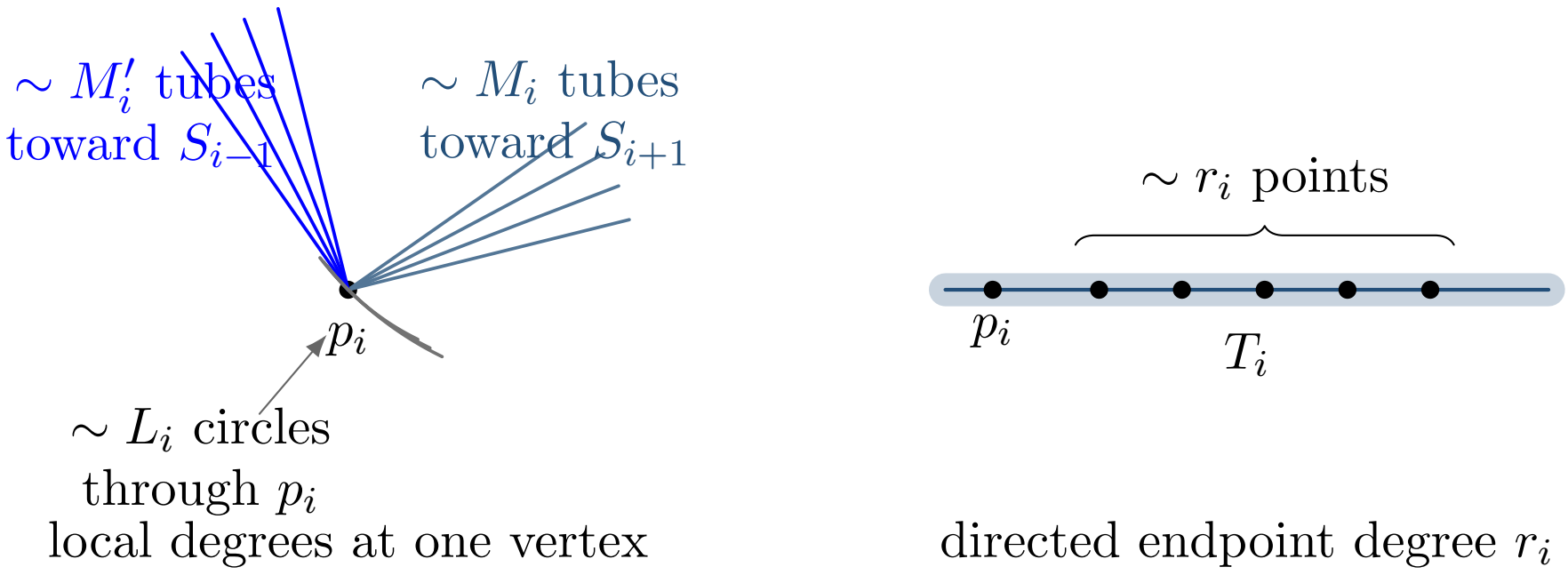
\begin{figure}[htbp]
		\centering
		\begin{tikzpicture}[scale=.88]
			\begin{scope}[xshift=-3.8cm]
				\node[point] (p) at (0,0) {};
				\node[lab,below] at (0,-.12) {$p_i$};
				\foreach \a in {14,21,28,35}{\draw[inkblue!78,line width=.7pt] (0,0)--({2.45*cos(\a)},{2.45*sin(\a)});}
				\foreach \a in {104,111,118,125}{\draw[blue,line width=.7pt] (0,0)--({2.45*cos(\a)},{2.45*sin(\a)});}
				\foreach \Rr in {1.55,1.82,2.09}{
					\draw[black!55,line width=.7pt] plot[domain=218:244,samples=35]
					({\Rr+1.414*\Rr*cos(\x)},{\Rr+1.414*\Rr*sin(\x)});
				}
				\node[lab,inkblue,align=left] at (1.75,1.48) {$\sim M_i$ tubes\\toward $S_{i+1}$};
				\node[lab,blue,align=right] at (-1.75,1.48) {$\sim M'_i$ tubes\\toward $S_{i-1}$};
				\draw[->,black!60] (-.75,-1.05)--(-.18,-.38);
				\node[lab,align=center] at (-1.25,-1.45) {$\sim L_i$ circles\\through $p_i$};
				\node[boxlabel] at (0,-2.15) {local degrees at one vertex};
			\end{scope}
			\begin{scope}[xshift=3.8cm]
				\draw[inkblue!25,line width=7pt] (-2.55,0)--(2.55,0);
				\draw[inkblue,line width=.8pt] (-2.55,0)--(2.55,0);
				\node[point] at (-2.15,0) {}; \node[lab,below] at (-2.15,-.05) {$p_i$};
				\foreach \x in {-1.25,-.55,.15,.85,1.55}{\node[point] at (\x,0) {};}
				\draw[decorate,decoration={brace,amplitude=4pt}] (-1.45,.38)--(1.75,.38)
				node[midway,above=5pt,lab] {$\sim r_{i}$ points};
				\node[lab] at (0,-.55) {$T_{i}$};
				\node[boxlabel] at (0,-2.15) {directed endpoint degree $r_i$};
			\end{scope}
		\end{tikzpicture}
		\caption{Multiplicities after the hypergraph refinement.\label{fig:hyp-local-parameters}}
	\end{figure}
    
    \step{Selection of a `good' corner}
    We upper-bound the number of rectangles using the parameters pigeonholed above. Starting with $S_i$ we may (attempt to) create a rectangle in a clockwise fashion as follows: 
    \begin{enumerate}
        \item Select $p \in S_i;$ there are $|S_i|$ choices.
        \item Select a tube joining $S_i$ to $S_{i+1};$ there are $\sim M_{i}$ choices.
        \item Select a point in $S_{i+1}$ lying on this tube; there are $\sim r_{i}$ choices.
        \item Select a tube through $p$ (if it exists) which makes an angle of $\pi/2 + O(\de^{1-O(\e)})$ with the tube selected in (2) lying in $\cT_{i-1} $; at most $\dem$ choices.
        \item Select a point of $S_{i-1}$ lying on this tube; $\sim r'_{i-1}$ choices. 
    \end{enumerate}
\bluecorrection{The remaining vertex is determined within $O(\de)$ by \eqref{eq.para}, so has only $O(1)$ possible choices.}
In total we have found at most 
\begin{equation}\label{eq.count}
    \lesssim \dem |S_i|M_{i}r_{i}r'_{i-1},
\end{equation}
rectangles. By double counting the incidences between $S_{i+1}$ and $\T_{i}$ we have
\begin{equation}
    |\T_{i}|r_{i} \lesssim \dem M'_{i+1}|S_{i+1}|.
\end{equation}
Likewise
$$|\T_{i-1}|r'_{i-1} \lesssim \dem M_{i-1}|S_{i-1}|$$
and so 
\eqref{eq.count} becomes
\begin{equation}
    \lesssim \dem |S_{i-1}||S_i||S_{i+1}|\frac{M_{i}M_{i-1}M'_{i+1}}{|\T_{i}||\T_{i-1}|}.
\end{equation}
Similarly, by counting anticlockwise we may obtain
\begin{equation}
    \lesssim \dem |S_{i-1}||S_i||S_{i+1}|\frac{M'_{i}M_{i-1}M'_{i+1}}{|\T_{i}||\T_{i-1}|},
\end{equation}
and so 
$$|\cR|^2 \lesssim \dem (|S_{i-1}||S_i||S_{i+1}|)^2
\frac{M_iM'_i M_{i-1}^2(M'_{i+1})^2}{|\T_{i}|^2|\T_{i-1}|^2}.$$
Therefore
$$|\cR|^8 \lesssim \dem\left(\prod_{i=1}^4|S_i|\right)^6
\prod_{i=1}^4 \frac{M_i^3 (M'_i)^3}{|\T_i|^2|\T_{i-1}|^2}.$$
Using $|S_i|\leq N$ and the comparison between $M_i$ and $M'_i$, we may therefore find $j = 1,2,3,4$ so that
\begin{equation}\label{eq.firstrectanglecount}
    |\cR| \lesssim \dem |S_j|N^2\frac{M^3_{j}}{|\T_{j}||\T_{j-1}|}.
\end{equation}
Without loss of generality, suppose $j=1$ and for ease of notation write
$$S = S_1;  \qquad \T = \T_{1}; \qquad \T' = \T_{4}; $$
$$ M= M_{1}; \qquad L = L_1; \qquad t = t_1.$$
\step{Applying Proposition \ref{prop.circle}}
Recall there are $L$ circles which pass through each $p \in S$ which intersect $\sim t $ points from $S.$ Double counting these incidences we have
$$ |\cC|t \lesssim \dem |S|L.$$
We note another way to bound $|\cR|.$ We attempt to find a rectangle as follows: fix a point in $S_2,$ fix a point in $S_4;$ at most $N^2$ choices. By our encoding this determines at most $O(1)$ circles $C \in \cC.$ Now we have $\sim t$ choices for the point in $S.$ Therefore
\begin{equation}\label{eq.one}
    |\cR| \lesssim \dem N^2t \lesssim \dem\frac{|S|N^2L}{|\cC|}.
\end{equation}
Another way to bound $|\cR|$ is as follows: fix one point from $S$ and one from $S_4;$ at most $|S|N$ choices. The possible perpendicular tubes give at most $\dem r_1$ choices for the point in $S_2$. Since a point of $S_2$ belongs to at most $\dem$ tubes through a fixed $p\in S$, double counting gives $Mr_1\lesssim\dem N$. We arrive at
\begin{equation}\label{eq.two}
    |\cR| \lesssim \dem |S|Nr_1 \lesssim \dem |S|N^2/M.
\end{equation}
Recalling \eqref{eq.firstrectanglecount}:
\begin{equation}\label{eq.three}
    |\cR| \lesssim \dem\frac{|S|N^2M^3}{|\T||\T'|}.
\end{equation}

We now apply Proposition \ref{prop.circle} and our three bounds above. The fixed arcs were selected in Step 1, and Step 2 supplies at least $\dep M$ tubes from each family and $\dep L$ circles at every point of $S$. Applying the proposition with these lower bounds costs only factors $\dem$. Recall that
   $$\Sigma=\frac{|\T||\T'|}{M^2}+\frac{(|\T|+|\T'|)|\cC|}{ML}.$$
Before our application, first record that
\begin{equation}
    \min \{|\T|,|\T'|\}^2 \leq |\T||\T'|,
\end{equation}
and
\begin{equation}
    \min\{|\T|,|\T'|\}^2(|\T| + |\T'|) \lesssim |\T|^{3/2}|\T'|^{3/2}.
\end{equation}
   An application of Proposition \ref{prop.circle} thus returns 
   \begin{equation}\label{eq.application}
       \de^sM^{2+\sigma}|S|^3
    \lesssim_\e\dem\Sigma\min\{|\T|,|\T'|\}^2 \lesssim \dem\left(\frac{|\T|^2|\T'|^2}{M^2}
    +\frac{(|\T||\T'|)^{3/2}|\cC|}{ML}\right).
   \end{equation}
We use \eqref{eq.three} to obtain
\begin{equation}
    \frac{|\T|^2|\T'|^2}{M^2} \lesssim \dem\frac{|S|^2N^4M^4}{|\cR|^2}.
\end{equation}
Using both \eqref{eq.one} and \eqref{eq.three} we bound
\begin{equation}
    \frac{(|\T||\T'|)^{3/2}|\cC|}{ML} \lesssim \dem\frac{|S|^{5/2}N^5M^{7/2}}{|\cR|^{5/2}}.
\end{equation}
Using the two bounds above with \eqref{eq.application}, along with \eqref{eq.two} and $0<\sigma\leq1$, we obtain
\begin{align*}
    \de^s
    &\lesssim_\e\dem\left(\frac{N^4M^{2-\sigma}}{|S||\cR|^2}
    +\frac{N^5M^{3/2-\sigma}}{|S|^{1/2}|\cR|^{5/2}}\right)\\
    &\lesssim_\e\dem\frac{|S|^{1-\sigma}N^{8-2\sigma}}{|\cR|^{4-\sigma}}\\
    &\lesssim_\e\dem\frac{N^{9-3\sigma}}{|\cR|^{4-\sigma}},
\end{align*}
where the last inequality uses $|S|\leq N$. \bluecorrection{Rearranging gives the claimed estimate in Proposition \ref{prop.rect}.}
\end{steps}
\begin{remark}
The original pre-AI approach by the authors used a diagonal tube as an entry of a hyperedge instead of a circle. The problem is that a side tube does not determine a diagonal tube up to bounded multiplicity, since the eccentricity of rectangles is variable.
\end{remark}

\section{Proof of Proposition \ref{prop.proj}}\label{sect.proj}

 Write $$F=\max_{\theta\in\Theta}\cn{\pi_\theta(P)}.$$ We immediately note, as before, that $\Theta$ is a $(\de,\sigma,O(|\Theta|^{1-\sigma}))$-KT set. Apply Lemma \ref{lem.twoends} to $\Theta.$ If the two-ends scale $\rho_0\geq\de^\e/100$, the rescaling below costs only factors $\dem$, and Theorem \ref{theorem.twoends} gives
$$F\gtrsim_{\e}\dep |\Theta|^{\sigma/2}\de^{s/2}|P|.$$
A geometric average with the bilinear bound $F\gtrsim\dep|P|^{1/2}$ gives \eqref{eq.propconcl}.

From now on we assume that $\rho_0 < \de^{\e}/100.$ Let $I$ be the interval given to us by Lemma \ref{lem.twoends}. \correction{By the assumption on $\Th$, we may choose $x_1,x_2\in\Theta$ with $|x_1-x_2|\gtrsim\depp$ and $d(x_j,I)\gtrsim\depp$ for $j=1,2$. Keep the original coordinates and write $\pi_x=\pi_{x_1}$, $\pi_y=\pi_{x_2}$. For $j\in\{x,y\}$, let $B_j(A)$ be the number of intervals in the fixed grid $\{[k\de,(k+1)\de):k\in\Z\}$ meeting $\pi_j(A)$. Thus $B_j(A)\sim\cn{\pi_j(A)}$, and $B_j(A)\lesssim F$ for $A\subset P$.}

\correction{We may assume that $\de/\rho_0$ is sufficiently small in terms of $\e$: otherwise $|\Theta\cap I|\lesssim_\e1$ by $\de$-separation, hence $|\Theta|\lesssim_\e\de^{-\e}$, and \eqref{eq.propconcl} follows from the bilinear bound and $|P|\lesssim\de^{-s}$.}

\correction{We require the following version of Lemma \ref{lem.rectangularpruninglol}.}
\begin{lemma}\label{lem.pruning}
There exists $P_1 \subset P$ with $|P_1| \geq \dep|P|$ so that the following holds for all $\de \leq \rho \leq 1.$ 

Let $Q \in \dy \rho {P_1}.$ Then 
\begin{equation}
     \correction{B_x(P \cap Q)B_y(P \cap Q)} \leq \dem \frac{F^2}{\cns \rho {P_1}}.
\end{equation}
\end{lemma}
\correction{\begin{proof}
Take a uniform subset $P_0\subset P$ with $|P_0|\gtrsim\de^\e|P|$. At every scale $\rho\geq\de$, a fixed pair of \bluecorrection{grid intervals} can occur in at most $O(\de^{-2\e})$ dyadic $\rho$-squares, since the $O(\rho)$-neighbourhoods of their inverse-image strips intersect inside a ball of radius $O(\rho\de^{-\e})$. Hence
\[
\sum_{Q\in\dy\rho{P_0}}B_x(P\cap Q)B_y(P\cap Q)\lesssim\de^{-2\e}F^2.
\]
Apply the argument of Lemma \ref{lem.rectangularpruninglol} to $P_0$, keeping these products computed using the original $P$. The surviving set $P_1$ has size $\gtrsim\dep|P|$ and satisfies the stated bound with denominator $N_\rho(P_1)$.
\end{proof}}
From now on replace and relabel $P$ with $P_1.$ We require a similar version of Proposition \ref{prop.RKT}. Many steps are identical and so we opt for a terser proof.
\begin{lemma}\label{lem.rectproj}
    Let $R$ be a rectangle pointing perpendicular to $I,$ of dimensions $\de \leq r \leq w \leq1,$ with $r \leq \rho_0.$ Then
    \begin{equation}
        |P \cap R| \leq C(\sqrt{rw}/\de)^s,
    \end{equation}
    where $$C \lesssim \dem \frac{F}{\sqrt{|P|}}.$$
\end{lemma}
\begin{proof}
    Let $\de \leq r \leq w,$ with $r \leq \rho_0,$ and let $R$ be a rectangle pointing perpendicular to $I$, with shorter side $r$ and longer side $w.$ The aim is to count
    $$|P \cap R|.$$
    (Cover $R$ by $O(1)$ dyadic $w$-squares; the argument below is applied in each such square $Q$, and the resulting bounds are summed.)
    We may cover \correction{the occupied part of} $R$ by $N \lesssim w/r$ $r$-squares from $\dy r P;$ call them $S_i$ for $1 \leq i \leq N.$ We then have, using the $(\de,s)$-KT condition on $P$ that
    \begin{align}
        |P \cap R| \lesssim (r/\de)^{s/2}  \sum_{i=1}^N |P \cap S_i|^{1/2}.
    \end{align}
    The aim is to now bound 
    \begin{equation}\label{eq.boundthisproj}
    \sum_{i=1}^N |P \cap S_i|^{1/2}.
        \end{equation}

 \correction{The fibres of $\pi_x$ and $\pi_y$ make an angle $\gtrsim\depp$ from the centre-line of $R$. We have}
 \begin{equation}
     \correction{|P \cap S_i|\lesssim\de^{-2\e}B_x(P \cap S_i)B_y(P \cap S_i).}
 \end{equation}
 \correction{Indeed, the inverse image of a pair of \bluecorrection{grid intervals} has diameter $O(\de^{1-\e})$ and contains at most $O(\de^{-2\e})$ points of the $\de$-separated set $P$. Also, since $\de\leq r$, the inverse image of each \bluecorrection{grid interval} meets only $O(\demm)$ of the squares $S_i$. We may therefore bound this sum by}
\begin{align}
     \sum_{i}|P\cap S_i|^{1/2}
     &\lesssim \dem \sum_{i}\correction{B_x(P \cap S_i)^{1/2}B_y(P \cap S_i)^{1/2}}\\
     &\leq \dem\sqrt{\correction{\sum_{i}B_x(P \cap S_i)\sum_{i}B_y(P \cap S_i)}}\\
     &\leq \dem \sqrt{\correction{B_x(P \cap Q)B_y(P \cap Q)}}\\
     &\leq \frac{\dem}{\sqrt {\cns w P}}\correction{F}.
\end{align}
Here we have used Cauchy--Schwarz and Lemma \ref{lem.pruning}. Since $P$ is $(\de,s)$-KT it follows that
\begin{equation}
    \cns w P \geq \dep \de^s |P| w^{-s}.
\end{equation}
We may further bound \eqref{eq.boundthisproj} by
\begin{equation}
    \dem (w/\de)^{s/2}\frac{F}{\sqrt{|P|}}.
\end{equation}
We may therefore finally bound
\begin{equation}
    |P \cap R| \leq C (\sqrt{rw}/\de)^s
\end{equation}
 where
 $$C \lesssim \dem \frac{F}{\sqrt{|P|}}.$$

\end{proof}
We are now in a position to complete the proof of Proposition \ref{prop.proj}. We follow the argument presented in the proof of \cite[Theorem 1.5]{doprod}.

\begin{steps}
    \step{Formulating the problem as an incidence problem.} \correction{For each $\theta\in\Theta\cap I$, take the inverse images of the occupied \bluecorrection{grid intervals of length $\delta$}. Let $\T$ be the resulting tube family, and let $\T(p)$ contain the one selected tube through $p$ in each direction.} Write $N=|\Theta\cap I|\mathrel{\correction{\gtrsim}}\depp |\Theta|$. Then $|\T|\lesssim NF$, and each $p\in P$ meets $\sim N$ tubes of $\T$.

    \correction{Tile the working region by half-open $\rho_0\times1$ rectangles $R^0$ pointing perpendicular to $I$, and let $\cR$ consist of sufficiently large fixed dilates $R$ of these cells. Assign points to their unique base cells and write $P_R$ for the assigned points. Put $\T_R=\bigcup_{p\in P_R}\T(p)$. All selected tubes lie in the corresponding $R$, and, since their directions lie in an arc of length $\rho_0$, each fine-tube label occurs in only $O(1)$ families $\T_R$. Thus}
    $$\sum_{R\in\cR}|\T_R|\lesssim|\T|,\qquad \sum_{R\in\cR}|P_R|\mathrel{\correction{=}}|P|.$$

    \step{Segments $\gamma$ and their multiplicity}
    Tile $R$ with a collection $\Gamma_{R,\mathrm{all}}$ of segments $\gamma$ with dimensions $(\delta,\delta/\rho_0)$, having the same orientation and eccentricity as $R$. \correction{Use fixed enlargements of the fine tubes in width and length when interpreting containment of these segments.}

    Each $p\in P_R$ lies inside a unique $\gamma\subset R$. We pick a dyadic $H$ such that, writing
    $$\Gamma_R=\{\gamma\in\Gamma_{R,\mathrm{all}}:|P_R\cap\gamma|\sim H\},$$
    we have
    $$\sum_{R\in\cR}|\Gamma_R|\gtrsim\dep\frac{|P|}{H}.$$
    We may recast the incidence graph as a new graph with vertex set $\T\sqcup\bigcup_R\Gamma_R$, and an edge between $\gamma$ and $T$ if and only if $\gamma\subset T$. \correction{If $p\in P_R\cap\gamma$, every selected $T\in\T(p)$ contains $\gamma$ after these fixed enlargements: the angular variation is $O(\rho_0)$ over length $O(\delta/\rho_0)$, costing only $O(\delta)$ extra width. Conversely, each direction contributes only $O(1)$ tubes containing $\gamma$. Thus each retained $\gamma$ has $\sim N$ incident tubes in $\T_R$, with all directions in $\Theta\cap I$ represented.}

    \step{Rescaling $R$}
    Let $\bar\delta=\delta/\rho_0$. Fix $R$. \correction{After rotation and translation, rescale $R$ onto $[0,1]^2$ with scaling factors comparable to $(1/\rho_0,1)$; all fixed normalisations only change absolute constants.} Then each $T\in\T_R$ becomes a $\bar\delta$-tube, and we call this collection $\bar\T_R$. Also, each $\gamma\in\Gamma_R$ becomes a $\bar\delta$-point, and we call this collection $\bar P_R$. Then $|\bar P_R|=|\Gamma_R|$ and $|\bar\T_R|=|\T_R|$. Our rescaling magnifies angles between tubes by a factor of $1/\rho_0$.

    \step{Preservation of spacing properties via rescaling}
    First, $\bar P_R$ is a $(\bar\delta,s,K_1)$-KT set with $K_1\lesssim C/(H\rho_0^{s/2})$, where $C$ is the constant in Lemma \ref{lem.rectproj}. Indeed, after rescaling, the number of $\bar p\in\bar P_R$ lying inside an $r$-ball with $r\geq\bar\delta$ is dominated by $1/H$ times the number of $p\in P$ that lie inside a rectangle with dimensions $(r\rho_0,r)$. By Lemma \ref{lem.rectproj}, this number is at most
    $$\frac{C}{H}\left(\frac{r\sqrt{\rho_0}}{\delta}\right)^s=\frac{C}{H\rho_0^{s/2}}\left(\frac{r}{\bar\delta}\right)^s.$$

    Second, for each $\bar p\in\bar P_R$, the family $\bar Y(\bar p)$ of incident tubes in $\bar\T_R$ is a $(\bar\delta,\sigma,K_2)$-KT set with $K_2\lesssim |\Theta|^{1-\sigma}$. Moreover, $|\bar Y(\bar p)|\sim N$. Its directions are the rescaled directions in $\Theta\cap I$, so Lemma \ref{lem.twoends} also gives \correction{the $(\bar\delta,\e,O(1))$-set condition, up to a uniformly bi-Lipschitz reparametrisation. Since $\bar\delta$ is sufficiently small in terms of $\e$, this implies the $(\bar\delta,\e,\bar\delta^{-\e^3})$-set condition required by Theorem \ref{theorem.twoends}}.

    \step{Applying the two-ends theorem}
    We apply Theorem \ref{theorem.twoends} to obtain
    \begin{align*}
        |\T_R|=|\bar\T_R|
        &\gtrsim_{\e}\dep K_1^{-1}K_2^{-1/2}N^{3/2}\bar\delta^{s/2}|\bar P_R|\\
        &\gtrsim_{\e}\dep\frac{H}{C}|\Theta|^{-(1-\sigma)/2}N^{3/2}\de^{s/2}|\Gamma_R|.
    \end{align*}
    Summing in $R$ gives
    $$NF\gtrsim_{\e}\dep C^{-1}|\Theta|^{-(1-\sigma)/2}N^{3/2}\de^{s/2}|P|.$$
    Since $N\mathrel{\correction{\gtrsim}}\depp |\Theta|$ and $C\lesssim\dem F/|P|^{1/2}$, we obtain
    $$F^2\gtrsim_{\e}\dep |\Theta|^{\sigma/2}\de^{s/2}|P|^{3/2},$$
    which proves \eqref{eq.propconcl} and completes the proof.
\end{steps}

\section{Fourth energy: conjectures, examples and decoupling}\label{sect.fourth-final}
Write $\Phi(x,y)=(x,y,x^2+y^2)$, and abbreviate
$\en{4}{\de}(P)=\en{4}{\de}(P,P,P,P)$. We present two conjectures and a sharp result regarding the number of approximate additive 4-tuples on the paraboloid.

\subsection{Two conjectures}
\begin{conj}[Separated fourth energy]\label{conj.fourth-separated}
Let $0<s\leq2$, and let $P_1,P_2,P_3,P_4\subset\bbP^2$ be $(\de,s)$-KT sets, each with cardinality comparable to $N=\de^{-s}$, mutually separated by a fixed distance $c>0$. For every $\eta>0$,
\begin{equation}\label{eq.conj-fourth-separated}
\en{4}{\de}(P_1,P_2,P_3,P_4)
\lesssim_{\eta,c}\de^{-\eta}(N^2+\de N^3)
\sim\de^{-a(s)-\eta},
\qquad a(s)=\max\{2s,3s-1\}.
\end{equation}
\end{conj}
In particular, the proposed exponent is $2s$ for $s\leq1$, and $3s-1$ for $s\geq1$. At $s=1$ the present result gives $7/3$, while the conjecture asks for $2$. The ``easy'' exponent $5/2$ is a consequence of Proposition \ref{prop.linear-fourth-decoupling} and H\"older's inequality.

\begin{conj}[Linear fourth energy]\label{conj.fourth-linear}
Let $0<s\leq2$ and set
\begin{equation}\label{eq.linear-fourth-prediction}
b(s)=\begin{cases}
5s/2,&0<s\leq1,\\
1+3s/2,&1\leq s\leq4/3,\\
3s-1,&4/3\leq s\leq2.
\end{cases}
\end{equation}
A $(\de,s)$-KT set $P\subset\bbP^2$ with $|P|\sim\de^{-s}$ should satisfy, for every $\eta>0$,
\begin{equation}\label{eq.conj-fourth-linear}
\en{4}{\de}(P)\lesssim_\eta\de^{-b(s)-\eta}.
\end{equation}
\end{conj}
The range $s\leq1$ is a theorem proved below, and a train-track example will show its sharpness. In this range, the exponent $\frac{10-3s}{4-s}$ of $\de^{-s}$ in the separated bound \eqref{eq.energyresult} is smaller than the linear exponent $5/2$.

Conjecture \ref{conj.fourth-linear} (and thus, Conjecture \ref{conj.fourth-separated}, too) also admits an easy resolution when $3/2\le s\le 2$, using the incidence result in  \cite{fu2021incidence}. Since the argument is rather standard, we do not include it here. We do not have a proof for the conjectured estimate in the intermediate range $1<s<3/2$.

\subsection{Examples}We discuss two examples for sharpness. 

\noindent\textbf{Arithmetic grids.}
Let $m\sim\de^{-s/2}$, $Q=m^{-1}\{0,\ldots,m-1\}^2$, and $P=\Phi(Q)$. Then $|P|=m^2\sim\de^{-s}$ and $P$ is $(\de,s)$-KT. Axis-parallel rectangles give $\gtrsim m^4$ quadruples with mutually separated vertices. For $s>1$, approximate rectangles give the further lower bound
\begin{equation}\label{eq.grid-fourth-lower}
\en{4}{\de}(P)\gtrsim m^4+\de m^6\sim\de^{-2s}+\de^{-(3s-1)}.
\end{equation}
Here is a direct count for the second term. Write an integer parallelogram as $p,p+u,p+u+v,p+v$. Take $\gtrsim m^2$ primitive integer vectors $u$ of length a small fixed multiple of $m$. For each integer $|k|\leq c_1\de m^2$, the line $u\cdot v=k$ has an integer solution $v$ in a fixed annulus of radius comparable to $m$, since its lattice spacing is $|u|$. Choose the annulus wide enough relative to $|u|$. There are also $\gtrsim m^2$ admissible translates $p$. All four vertices can be kept in the grid and mutually separated by a fixed multiple of $m$. This gives $\gtrsim\de m^6$ quadruples. Partitioning the bounded patch into finitely many sufficiently small camps retains a fixed fraction in four separated camps. Thus both branches of \eqref{eq.conj-fourth-separated} are necessary; see Figure \ref{fig:fourth-grid}.

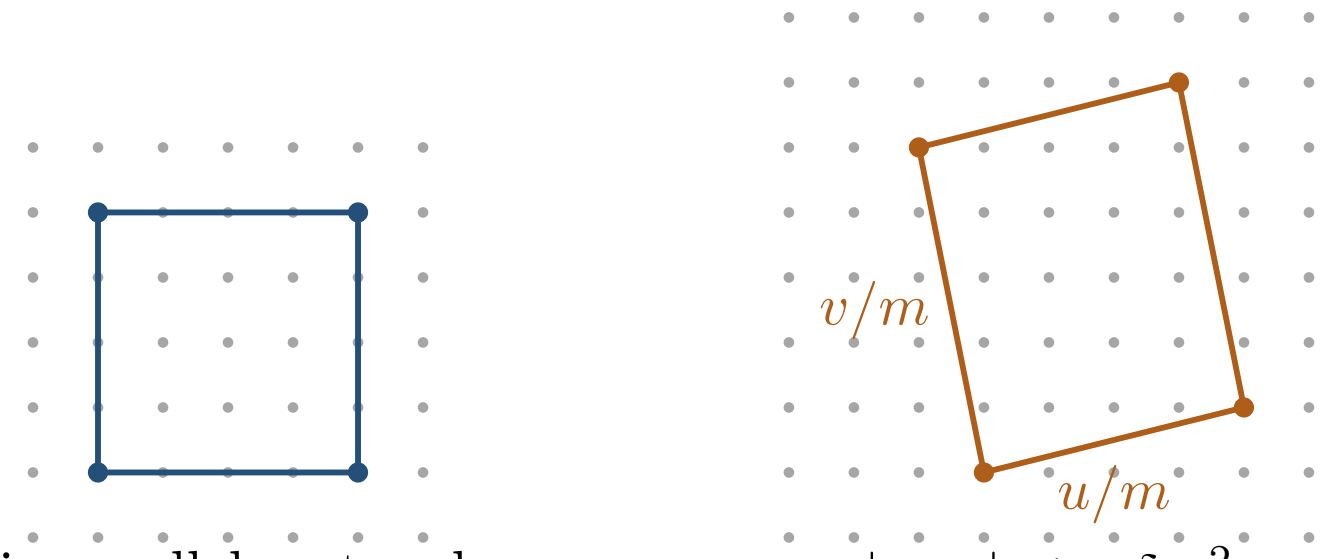
\begin{figure}[htbp]
\centering
\begin{tikzpicture}[x=.43cm,y=.43cm]
  \begin{scope}
    \foreach \x in {0,...,6}{\foreach \y in {0,...,6}{
      \fill[black!35] (\x,\y) circle (1pt);}}
    \draw[inkblue,line width=1.1pt] (1,1)--(5,1)--(5,5)--(1,5)--cycle;
    \foreach \p in {(1,1),(5,1),(5,5),(1,5)}{\fill[inkblue] \p circle (1.9pt);}
    \node[lab,align=center] at (3,-1.15) {axis-parallel rectangles\\$\gtrsim m^4$ quadruples};
  \end{scope}
  \begin{scope}[xshift=5cm]
    \foreach \x in {0,...,8}{\foreach \y in {0,...,8}{
      \fill[black!35] (\x,\y) circle (1pt);}}
    \draw[inkorange,line width=1.1pt] (3,1)--(7,2)--(6,7)--(2,6)--cycle;
    \foreach \p in {(3,1),(7,2),(6,7),(2,6)}{\fill[inkorange] \p circle (1.9pt);}
    \node[lab,inkorange,below] at (5,1.45) {$u/m$};
    \node[lab,inkorange,left] at (2.5,3.5) {$v/m$};
    \node[lab,align=center] at (4,-1.15) {$|u\cdot v|\leq c_1\de m^2$\\$\gtrsim\de m^6$ quadruples};
  \end{scope}
\end{tikzpicture}
\caption{Grid examples in the parameter plane, drawn schematically. The four lifts of $(p,p+u,p+u+v,p+v)/m$ have additive error $2|u\cdot v|/m^2$. Exact rectangles give the first term in \eqref{eq.grid-fourth-lower}; near-right-angle parallelograms give the second.}
\label{fig:fourth-grid}
\end{figure}

\medskip\noindent\textbf{Train tracks.}
Put
\[
t=\min\{s,1\},\quad \rho=\sqrt\de,\quad
m\sim\de^{-t/2},\quad n\sim\de^{-s/2},\quad K\sim\de^{-(s-t)/2},
\]
and, with a sufficiently small fixed $c_2>0$, take
\[
A=\left\{\frac{k}{2K}+\frac{c_2\rho i}{m}:0\leq k<K,\ 0\leq i<m\right\},
\quad B=\{j/n:0\leq j<n\}.
\]
Set $P=\Phi(A\times B)$. Choose fixed comparison constants so that $P$ is $\de$-separated. Integer parts do not affect any powers. There are $K$ horizontal clusters, each of width $c_2\rho$, and $|P|=Kmn\sim\de^{-s}$. Counting points in balls verifies the KT condition. For $s>1$, the relevant scales are $\de\leq n^{-1}\leq\rho\leq K^{-1}$; on the successive ranges the planar ball counts are bounded by
\[
C r/\de,\quad C r^2\de^{-1-s/2},\quad
C r\de^{-(s+1)/2},\quad C r^2\de^{-s},
\]
respectively, and each is at most $C(r/\de)^s$. For $s\leq1$, $K=1$; below $\rho$ a ball meets $O(1)$ rows and above $\rho$ it meets $O(1+rn)$ rows, giving the same KT conclusion.

Choose all four horizontal coordinates in one cluster, with $i_1-i_2+i_3-i_4=0$, and choose vertical indices with $j_1=j_2$ and $j_3=j_4$, as in Figure \ref{fig:fourth-tracks}. Horizontal closure is exact; the remaining quadratic error is $O(c_2^2\de)$. Consequently
\begin{equation}\label{eq.train-fourth-lower}
\en{4}{\de}(P)\gtrsim Km^3n^2
\sim\de^{-(3s/2+\min\{s,1\})}.
\end{equation}
Together, \eqref{eq.grid-fourth-lower} and \eqref{eq.train-fourth-lower} force exactly $b(s)$. These train-track quadruples have a pair at distance $O(\sqrt\de)$, so they do not obstruct Conjecture \ref{conj.fourth-separated}.

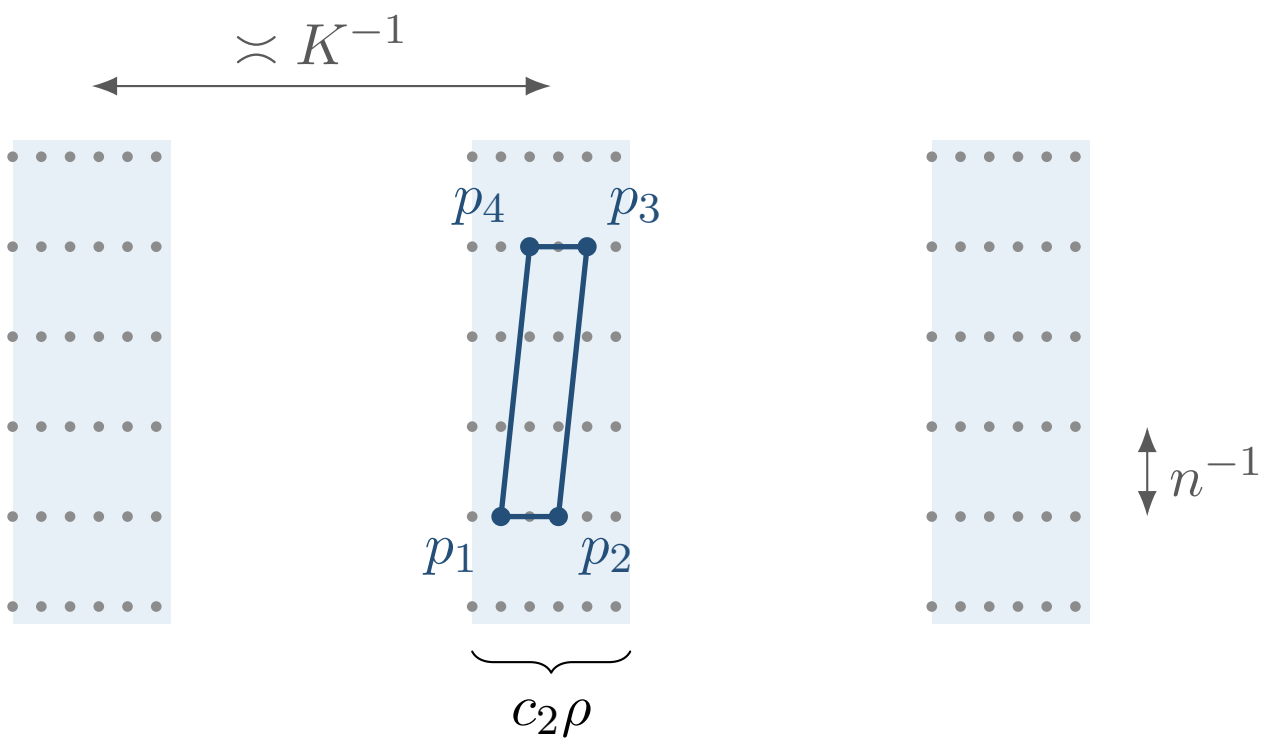
\begin{figure}[htbp]
\centering
\begin{tikzpicture}[x=.95cm,y=.85cm]
  \foreach \a in {0,3.2,6.4}{
    \fill[softblue] (\a,-.13) rectangle ({\a+1.1},3.63);
    \foreach \i in {0,...,5}{\foreach \j in {0,...,5}{
      \fill[black!45] ({\a+.2*\i},{.7*\j}) circle (1pt);}}
  }
  \draw[inkblue,line width=1pt] (3.4,.7)--(3.8,.7)--(4,2.8)--(3.6,2.8)--cycle;
  \foreach \x/\y/\name/\where in {3.4/.7/$p_1$/below left,3.8/.7/$p_2$/below right,4/2.8/$p_3$/above right,3.6/2.8/$p_4$/above left}{
    \fill[inkblue] (\x,\y) circle (1.8pt);
    \node[lab,\where,inkblue] at (\x,\y) {\name};}
  \draw[decorate,decoration={brace,amplitude=4pt,mirror}] (3.2,-.35)--(4.3,-.35)
    node[midway,below=5pt,lab] {$c_2\rho$};
  \draw[<->,black!65] (.55,4.05)--(3.75,4.05)
    node[midway,above,lab] {$\asymp K^{-1}$};
  \draw[<->,black!65] (7.9,.7)--(7.9,1.4)
    node[midway,right,lab] {$n^{-1}$};
  \node[lab] at (3.75,-1.5) {$K$ bands, $m$ columns per band, $n$ rows};
\end{tikzpicture}
\caption{Train tracks in the parameter plane. The highlighted quadruple has $y_1=y_2$, $y_3=y_4$, and $x_1-x_2+x_3-x_4=0$. Its quadratic error after lifting is $O(c_2^2\de)$, since all four horizontal coordinates lie in one band of width $c_2\rho$. When $s\leq1$, there is only one band.}
\label{fig:fourth-tracks}
\end{figure}

\subsection{The decoupling bound}
The following immediate application of decoupling verifies Conjecture \ref{conj.fourth-linear} when $s\le 1$. This approach becomes inefficient for larger values of $s$.

\begin{proposition}\label{prop.linear-fourth-decoupling}
Let $0<s\leq2$ and $\eta>0$. Every $(\de,s)$-KT set $P\subset\bbP^2$ satisfies
\begin{equation}\label{eq.linear-fourth-decoupling}
\en{4}{\de}(P)\lesssim_\eta\de^{-s/2-\eta}|P|^2
\lesssim_\eta\de^{-5s/2-\eta}.
\end{equation}
\end{proposition}
\begin{proof}
Set $R=\de^{-1}$. Partition the paraboloid into caps $\theta$ of tangential diameter $R^{-1/2}$, and write
\[
S(\xi)=\sum_{p\in P}e^{-2\pi i\xi\cdot p},\qquad
S_\theta(\xi)=\sum_{p\in P\cap\theta}e^{-2\pi i\xi\cdot p},\qquad
n_\theta=|P\cap\theta|.
\]
The KT condition gives $n_\theta\lesssim R^{s/2}$ and $\sum_\theta n_\theta=|P|$.
Let $w_R(\xi)=(1+|\xi|^2/R^2)^{-100}$. Its Fourier kernel decays rapidly, so $\de$-separation gives
\[
\|S_\theta\|_{L^2(w_R)}^2
\lesssim R^3\sum_{p,q\in P\cap\theta}(1+R|p-q|)^{-10}
\lesssim R^3n_\theta.
\]
Since $\|S_\theta\|_\infty\leq n_\theta$,
$\|S_\theta\|_{L^4(w_R)}^2\lesssim R^{3/2}n_\theta^{3/2}$.
Critical $\ell^2L^4$ decoupling for the elliptic paraboloid \cite[Theorem 1.1]{BD}, in its weighted localised form, therefore yields
\[
\int w_R|S|^4
\lesssim_\eta R^{3+\eta}\left(\sum_\theta n_\theta^{3/2}\right)^2
\lesssim_\eta R^{3+s/2+\eta}|P|^2.
\]
Gaussian Fourier inversion gives
\[
\int e^{-\pi|\xi|^2/R^2}|S(\xi)|^4\,d\xi
=R^3\sum_{p_1,p_2,p_3,p_4\in P}
e^{-\pi R^2|p_1-p_2+p_3-p_4|^2}
\gtrsim R^3\en{4}{1/R}(P).
\]
The Gaussian is dominated by a fixed multiple of $w_R$, proving \eqref{eq.linear-fourth-decoupling}.
\end{proof}
For $s\leq1$, the train-track lower bound \eqref{eq.train-fourth-lower} matches Proposition \ref{prop.linear-fourth-decoupling} up to the factor $\de^{-\eta}$, proving Conjecture \ref{conj.fourth-linear} in this range.

\label{sect.fourth-final-end}

\makeatletter
\let\blueaudit@lbibitem\@lbibitem
\def\@lbibitem[#1]#2{%
  \def\blueaudit@key{#2}\def\blueaudit@orp{O1}\def\blueaudit@wz{wang2022sticky}%
  \ifx\blueaudit@key\blueaudit@orp\color{black}%
  \else\ifx\blueaudit@key\blueaudit@wz\color{black}\else\color{black}\fi\fi
  \blueaudit@lbibitem[#1]{#2}}
\makeatother

\enlargethispage{2\baselineskip}
    	\bibliographystyle{alpha}
	\bibliography{references}

    \bigskip
\noindent\textsc{Ciprian Demeter}\\
Department of Mathematics, Indiana University\\
831 East 3rd St., Bloomington, IN 47405\\
\textit{Email address:} \href{mailto:demeterc@iu.edu}{demeterc@iu.edu}

\medskip
\noindent\textsc{William O'Regan}\\
Department of Mathematics, The University of British Columbia\\
Room 121, 1984 Mathematics Road, Vancouver, BC, Canada V6T 1Z2\\
\textit{Email address:} \href{mailto:woregan@math.ubc.ca}{woregan@math.ubc.ca}

    \end{document}